\documentclass[a4paper,12pt]{article}
\usepackage{graphicx}
\usepackage{float}
\usepackage{wrapfig}
\usepackage{amsfonts}
\usepackage{eurosym}
\usepackage{amssymb}
\usepackage{amsmath}
\usepackage{subcaption} 
\usepackage{natbib}
\usepackage[a4paper]{geometry}
\usepackage[singlespacing]{setspace}
\newtheorem{theorem}{Theorem}

\newtheorem{definition}[theorem]{Definition}

\newtheorem{proposition}{Proposition}

\newenvironment{proof}[1][Proof]{\noindent\textbf{#1.} }{\ \rule{0.5em}{0.5em}}
\begin{document}

\title{Increasing Risk: Dynamic Mean-Preserving Spreads}

\author{Jean-Louis Arcand\thanks{{\small Centre for Finance and Development
and Department of International Economics, The Graduate Institute, Geneva,
Switzerland. \ Email: }\textsf{jean-louis.arcand@graduateinstitute.ch}} \and %
Max-Olivier Hongler\thanks{{\small EPFL-IPR-LPM, Ecole Polytechnique F\'{e}d%
\'{e}rale de Lausanne, Lausanne, Switzerland. \ Email: }\textsf{%
max.hongler@epfl.ch}} \and Daniele Rinaldo\thanks{%
Department of International Economics, The Graduate Institute, Geneva,
Switzerland. Email: daniele.rinaldo@graduateinstitute.ch}}
\date{May 3, 2017}
\maketitle

\begin{abstract}
We extend the celebrated Rothschild and Stiglitz (1970) definition of Mean-Preserving Spreads to a dynamic framework. We adapt the original integral conditions to transition probability densities, and give sufficient conditions for their satisfaction. We then prove that a specific nonlinear scalar diffusion process, super-diffusive ballistic noise, is the unique process that satisfies the integral conditions among a broad class of processes. This process can be generated by a random superposition of linear Markov processes with constant drifts. This exceptionally simple representation enables us to systematically revisit, by means of the properties of Dynamic Mean-Preserving Spreads, four workhorse economic models originally  based on White Gaussian Noise.     
\end{abstract}

\onehalfspacing

\newpage

\section{Introduction}

Comparing the riskiness of different random variables is a topic of central importance in economic research. The inadequacy of the variance as a measure of risk is well established, since this criterion is satisfactory in economic applications in a limited number of cases. To wit: an increase in risk increases the variance, but the converse is not necessarily true. \\ 

A milestone in the search for a more informative criterion was the series of articles by Rothschild and Stiglitz (\citeyear{Rothschild70}, \citeyear{Rothschild71}, \citeyear{Rothschild72}) which defined the concept of an increase in risk in the form of second-order stochastic dominance, often referred to as a Mean-Preserving Spread (MPS)%
, and applied it to various economic problems.\footnote{%
Various terminologies apply. \ For example, \cite{Chateauneuf2004}
refer, more accurately, to a mean-preserving increase in risk (MPIR).} This concept has become a workhorse of microeconomic analysis, with
applications ranging from finance to the study of inequality: see, for example, the standard textbooks by \cite{Laffont90}, \cite{Levy98} or \cite{Gollier2001}. \\ 

The strength of Rothschild and Stiglitz's result lies in a definition of comparative risk that can be summarized by means of four intuitively-appealing notions which are shown to be equivalent. 
\begin{definition}\textsc{Rothschild and Stiglitz definition of increasing risk.}  \\
For two random variables $x$ and $y$ with arbitrary distributions $F(.)$ and $G(.)$,  $y$ is said to be riskier than $x$ if: \\

1.1 The distribution $G(.)$ can be derived by adding zero-mean ``noise'' to $F(.)$;\\

1.2 The distribution $G(.)$ can be derived from $F(.)$ by means of one or more mean-preserving spreads, i.e. $G(.)$ has more weight in the tails;\\

1.3 The distributions $F(.)$ and $G(.)$ obey two integral conditions, one that imposes first-order stochastic equivalence and the second that defines second-order stochastic dominance;\\

1.4 Any optimizing rational agent with concave utility will prefer $F(.)$ to $G(.)$.\\
\label{def1}
\end{definition}
More explicitly, and for comparison purposes with what follows, we focus on Definition 1.3, which is given by the two following integral conditions:\footnote{Since in this paper we deal solely with diffusion processes, Definition 2 corresponds to the differentiable version of the integral conditions, as given by \cite{Diamond74}. }

\begin{definition} \textsc{MPS Integral Conditions:}\\
Let $x\in \mathbb{R}$ be a random variable distributed according to the $\lambda$-parameterized continuous CDF $x\sim \mathcal{P}^{(\lambda)} (x)$ in a probability space $(\Omega, \mathcal{F}, P)$, with $\Omega = \infty$. Then an increase in $\lambda$ generates a mean-preserving spread if: 
\begin{eqnarray}
\frac{\partial}{\partial \lambda} \int_\mathbb{R}  \mathcal{P}^{(\lambda)} (x) dx&=& 0,\label{MPSRS1} \\
\frac{\partial}{\partial \lambda} \int_{-\infty}^y  \mathcal{P}^{(\lambda)} (x) dx &\geq& 0 \quad \forall y < +\infty. \label{MPSRS2}
\end{eqnarray}
\label{def2}
\end{definition}
The support of $\mathcal{P}^{(\lambda)} $ may obviously also be compact, in which case the limits of the integrals and the upper bound of $y$ correspond to the boundaries of the support.  We present this form of the integral conditions for consistency with what follows.\footnote{The four notions of Definition 1 have been further expanded by \cite{machinapratt} who define more general rules in the creation of sequences of MPSs and of the zero-conditional mean noise, as well as generalizing mean-preserving spreads beyond distributions which are discrete or that possess a well-defined density function. } \\

That risk and variance do not necessarily coincide, and that  a risk-averse agent does not necessarily prefer a distribution with a lower variance to one with a higher variance (where the means are the same), is often forgotten in applied economic research. The ``risk" faced by peasants in developing countries is often proxied by the variance of their crop yields. The ``risk" faced by an investor is often proxied by the variance of asset returns.  But Definition 1.4 involves risk, not variance. Consider a slightly extended version of the simple example provided by \cite{Laffont90}, p. 26. There are two lotteries $x_1$ and $x_2$, given by:
\begin{eqnarray*}
x_1 &=& [(0.01, 0.10); (0.10,0.00);(1,0.70); (10,0.00);(100, 0.20);(1090,0.00)], \\
x_2 &=& [(0.01, 0.00);(0.10, 0.01); (1,0.00);(10,0.98); (100,0.00);(1090, 0.01)],
\end{eqnarray*}
where each pair $(x_{ij}, p_{ij})$ corresponds to the probability $p_{ij}$ of the realization $x_{ij}$, for lottery $j =1,2$ and states of nature $i = 1,2,3,4,5,6$. Notice that the expected values of the two lotteries are the same: $\mathbb{E} x_1 = \mathbb{E} x_2 = 20.701$. However, $Var \ x_1 = 2,000.7 < 11,979 = Var \ x_2$. Despite the variance of $x_2$ being much larger than the variance of $x_1$, an agent with logarithmic utility will strictly prefer $x_2$ over $x_1$ because  $\mathbb{E} \log x_1 = 0.46 < 2.303 = \mathbb{E} \log x_2$. The reason is clear: while the discrete version of the first integral condition \eqref{MPSRS1} is satisfied, the second integral condition (equation 2) is not. To see why, notice that if we ``stop" at $y=0.10$, the second integral condition reads $p_{11} +p_{21}= 0.10 > 0.01 = p_{12}+p_{22}$, while if we ``stop" at $y=10$ the condition reads $p_{11}+p_{21} +p_{31}+p_{41}= 0.71 < 0.99  =p_{12}+p_{22}+p_{32}+p_{42}$. As such, the two lotteries cannot unambiguously be ranked in terms of their risk, while they can be in terms of their variance.\\

%
Definition 2 applies to a static framework: loosely speaking, the comparison of riskiness of random variables is done for a ``snapshot'' of their respective distributions taken at an arbitrary  instant in time, as in a phase diagram for a dynamical system. In this paper, we provide the dynamic counterpart to mean-preserving spreads in the context
of scalar diffusion processes.  This allows us to parameterize the riskiness of a stochastic process throughout its evolution in the time domain. A remarkable feature of our dynamic counterpart is that it allows one to prove, for any process that exhibits the Brownian bridge property, that a specific
functional form, which corresponds to super-diffusive ballistic noise,
constitutes the \textit{sole} process with non-constant drift that displays
the dynamic version of the MPS property. In what follows we refer to this as a \emph{Dynamic Mean-Preserving Spread}, or DMPS. While the functional form is non-Gaussian,
its properties allow for simple closed-form solutions in a broad
range of economic applications, of which we give four canonical examples below. \\

This paper is organized as follows.     In Section 2, we derive our main
results.  First, in Definition 3, we give the two necessary integral conditions for a DMPS, which are straightforward dynamic generalizations of the standard Rothschild and Stiglitz conditions of Definition 2. The two conditions are essentially  antisymmetry and positivity conditions on the derivative with respect to a risk parameter of the Radon-Nikodym derivative associated with the transition probability density that defines a family of risk-parameterized scalar diffusion processes.  In Proposition 1, we provide a sufficient condition for a stochastic process to satisfy the integral conditions of Definition 3.   This is followed by Lemma 1, which shows that Definition 3 allows one to characterize second-order stochastic dominance in terms of the preferences of a risk-averse agent, as in Definition 1.4 above.  
\\

Our most important result is given in Proposition 2 which proves that, among processes that exhibit the Brownian bridge property, the diffusion process given by $ dX_t =  \sqrt{2\lambda} \tanh( \sqrt{2\lambda} X_t )dt + dW_t$, where $W_t$ is the standard Brownian motion, is the \textit{only} process with non-constant drift that displays the DMPS property.\footnote{The function $\tanh (x)$ is the hyperbolic tangent function given by $\frac{e^x - e^{-x}}{e^x+e^{-x}}.$} It turns out that this process has an extremely simple
representation in terms of the superposition of two drifted Wiener
processes: we prove this in Lemma 2, which we call the Bernoulli Representation Lemma.  This leads to particularly simple closed-form solutions in common
applications.  To give a first taste of this underlying simplicity, Proposition 3 then uses the preceding results to provide the marginal densities of DMPS-driven processes for three cases often used in economics: the drifted process with scalar coefficients, the mean-reverting (stationary Ornstein-Uhlenbeck) process and the geometric process.  Section 2 concludes with Proposition 4 in which we derive It\^o's formula for a DMPS process. \\ 

Section 3 explores how driving a system with the DMPS noise process, and increasing its parameter of risk $\lambda$, differs from an increase in the variance of the Brownian motion for a general diffusion process. We do this in three ways.  First, in Proposition 5, we study the curvature of the time-invariant probability measure for scalar processes, and show that the behavior obtained by driving the system with the DMPS process cannot be derived from a simple change in the variance of a Gaussian.  To wit: varying the risk parameter $\lambda$ and the variance term induce very different effects on the stationary probability measure.  Second,  we study two simple optimal stopping applications (stopping at the ultimate maximum and stopping with a transaction cost), and show that the impact of an increase in risk on both the stopping threshold and the stopping region is different from that of an increase in variance.  Third, we show that, contrary to
an increase in variance, a DMPS may violate the Certainty Equivalence
Principle used in optimal control theory.  The upshot of Section 3 (as with the rest of the paper) is that, as in the static world of Rothschild and Stiglitz, risk and variance should not be conflated in a dynamic context.\ 
\\

Section 4 provides economic illustrations of our results and showcases the analytical tractability of this class of processes by revisiting four standard economic problems: (i) portfolio selection, (ii) investment under uncertainty as in \cite%
{Abel1983} and \cite{Abel1994}, (iii) asset dynamics \`a la Black-Scholes and
finally (iv) firm entry and exit decisions under uncertainty following the \cite{Dixit1989} framework.  Our goal with these illustrations is not to propose new theoretical models, but to show how driving noise with a DMPS process instead of a Gaussian, and thereby disentangling risk and variance, modifies and often clarifies existing results.  The Gaussian setup always emerges as a special, and sometimes misleading, case.

\newpage 
\section{Main results}

On $\mathbb{R}$, consider the scalar diffusion process $X_{t}$ defined
by the stochastic differential equation on a filtered probability space $(\Omega, \mathcal{F}, P)$:%
\begin{equation}
\left\{ 
\begin{array}{l}
dX_{t}=b(X_{t})dt+\sigma dW_{t}, \\ 
X_{0}=x_{0},%
\end{array}%
\right.  \label{SDA}
\end{equation}%
where we assume $X_t \in \mathbb{R}, t \in [0, T]$, $\sigma \in \mathbb{R}$, the measurable function $b: \mathbb{R} \to \mathbb{R}$ is assumed at least $\mathcal{C}^2$ with bounded first and second derivatives, and $W_{t}$ is the standard Brownian Motion. Associated with equation
(\ref{SDA}), we define, for any function $\varphi:\mathbb{R}\rightarrow 
\mathbb{R}^{+},$ the diffusion operator $\mathcal{L}$:%
\begin{equation}
\mathcal{L}_{x}\varphi (x):=\left[ {\frac{\sigma ^{2}}{2}}{\frac{\partial
^{2}}{\partial x^{2}}}+b(x){\frac{\partial }{\partial x}}\right] \varphi (x).
\label{DOP}
\end{equation}%
For any time $0 \leq s < t < T$ let us write the transition probability density (TPD)
which describes the diffusion process of equation (\ref{SDA}), as $q(x,t | x_0, s)$. Assume
that $H(x,t)$ is a positive classical solution of the partial differential equation:%
\begin{equation}
{\frac{\partial }{\partial t}}H(x,t)+\mathcal{L}_{x}\left[ H(x,t)\right] =0.
\label{CONDO}
\end{equation}%
Then applying It\^{o}'s Lemma to $H(X_{t},t)$ with the process $X_{t}$
defined by (\ref{SDA}), we have $\mathbb{E}\left\{ {\frac{d}{dt}}%
H(X_{t},t)=0\right\} $, where $\mathbb{E} \{ \cdot  \} $ stands for the
expectation operator. Let $\hat{X}_t$ be a weak solution of \eqref{SDA}; then, by Theorem 2.1 of \cite{DaiPra1991}, the stochastic differential equation:%
\begin{equation}
\left\{ 
\begin{array}{l}
d\hat{X}_{t}=\left\{ b(\hat{X}_{t})+\sigma ^{2}{\frac{\partial }{\partial x}}%
\log \left[ H(x,t)\right] \mid _{x=\hat{X}_{t}}\right\} dt+\sigma dW_{t}, \\ 
\,\hat{X}_{0}=x_{0},%
\end{array}%
\right.  \label{DUALSDA}
\end{equation}%
admits a solution in $[0,T]$.
The TPD $Q(x,t | x_0, s)$ characterizing the diffusion process $\hat{X}_{t}$ given by \eqref{DUALSDA} reads:%

\begin{equation}
Q(x,t | x_0, s)=\left[ \frac{H(x,t)}{H(x_{0},s)]}\right] q(x,t | x_0, s).
\label{DUALTPD}
\end{equation}%
The function $z_{t}:=\left[ \frac{H(x,t)}{H(x_{0},s)]}\right] $ is the
Radon-Nikodym derivative for the change of measure in $\Omega$ relating the TPD $Q(x,t | x_0, s)$ with $%
q(x,t | x_0, s)$ and the process $Z_{t}:=\left[ \frac{H(X_{t},t)}{H(X_{0},s)]} \right ] $ is a martingale with $\mathbb{E}\left \{  Z_{t}\right \} =1$ (\citeauthor{DaiPra1991}, 1991) with respect to the natural filtration $\mathcal{F}_t = \sigma \{ X_s: 0 \leq s \leq t \}$.\\

Consider the class of positive definite functions:

\begin{equation}
H(x,t)=e^{-\lambda t}h^{(\lambda) }(x),\qquad h^{(\lambda )}(x)\geq 0,\qquad
x\in \mathbb{R},  \label{SUBCLASS}
\end{equation}%
which in view of equation (\ref{CONDO}) implies:%
\begin{equation}
\mathcal{L}_{x}\left[ h^{(\lambda) }(x)\right] =\lambda \,h^{(\lambda )}(x),
\label{STURM}
\end{equation}%
where $\lambda \in \mathbb{R}^+ $ is a positive constant, which will correspond in what
follows to the Rothschild and Stiglitz parameter of increasing risk. \
Substituting equation (\ref{SUBCLASS}) into equation (\ref{DUALTPD}), we can
write a $\lambda $-family of TPDs as:%

\begin{equation}
Q^{(\lambda )}(x,t | x_0, s)=e^{-\lambda t}\left[ {\frac{h^{(\lambda )}(x)}{%
h^{(\lambda)}(x_{0}))}}\right] q(x,t | x_0, s).  \label{DUALTPDA}
\end{equation}%
Since $Z_t$ is a martingale with $\mathbb{E}\left \{ Z_{t}\right \}  =1$, equation (%
\ref{DUALTPDA}) itself defines a normalized TPD. With $x_{0}=0$ and $s=0$,
the mean of $\hat{X}_{t}$ is given by:%

\begin{equation}
m^{(\lambda) }(t)=\mathbb{E}^{(\lambda )}\left\{ \hat{X}_{t}\right\} ={\frac{%
e^{-\lambda t}}{h^{(\lambda) }(0)}}\int_{\mathbb{R}}xh^{(\lambda)
}(x)q(x,t| 0,0)dx.  \label{MEANA}
\end{equation}%
Let us now assume that, in equation (\ref{SDA}), we have $b(x)=-b(-x)$. In
view of equation (\ref{DOP}), this implies symmetry: $\mathcal{L}_{x}\equiv 
\mathcal{L}_{-x}$. In turn, equation (\ref{STURM}) implies that $%
q(x,t| 0,0)=q(-x,t| 0,0)$ and $h^{(\lambda) }(x)=h^{(\lambda) }(-x),$ and
therefore the antisymmetry of the integrand in equation (\ref{MEANA}). \ It
follows that $m^{(\lambda )}(t)\equiv 0$ for all $\lambda \in \mathbb{R}^+$ i.e. the first moment is unchanged by a variation in $\lambda$. \\

Let us now present the main definition of the paper.  A Dynamic Mean-Preserving Spread (or dynamic mean-preserving increase in
risk) with respect to the parameter of increasing risk $\lambda $ is defined
by the dynamic counterparts of the two well-known integral conditions of Rothschild and Stiglitz given in \eqref{MPSRS1} and \eqref{MPSRS2}. 

\begin{definition}\textsc{DMPS Integral Conditions.}\\
Define the transition cumulative density (TCD) to be given by: 
$$
\mathcal{P}^{(\lambda
)}(x,t):=\int_{-\infty }^{x}Q^{(\lambda)}(y,t|0,0)dy.
$$
Then a Dynamic
Mean-Preserving Spread (DMPS) is defined by:%

\begin{eqnarray}
i) &\ &\frac{\partial }{\partial \lambda }\left[ \int_{\mathbb{R}}%
\mathcal{P}^{(\lambda )}(x,t)dx\right] =0,   \qquad (antisymmetry) \label{MPSCONDO}\\ 
ii) &\ & \frac{\partial }{\partial \lambda }\left[ \int_{-\infty }^{x}%
\mathcal{P}^{(\lambda) }(y,t)dy\right] \geq 0, \qquad  (positivity)%
\label{MPSCONDO2}
\end{eqnarray}
for all $x \in \mathbb{R}$, $\lambda \in \mathbb{R}^+$ and $t \in [0,T]$ .
\label{def3}
\end{definition}
\bigskip
The integral conditions in Definition \ref{def3} are essentially identical to the
integral conditions in Definition \ref{def2}  except that instead of
a cumulative density as the integrands we now have a transition cumulative density which evolves over time: we have therefore extended the ``static" result of Rothschild and Stiglitz to a dynamic framework.\\ 

Let us further explain this point. The original conditions of Definition 2 allowed one to parameterize the riskiness of different distributions by means of a  partial ordering in terms of second-order stochastic dominance: an increase in a single parameter implies an increase in the risk of the distribution. The two conditions of Definition 3 reflect the same goal, but allow for the ordering to be extended to time-evolving stochastic processes. With the generalization of the Rothschild-Stiglitz integral conditions to scalar diffusion processes, one can build a framework where risk and variance are disentangled in dynamic contexts. Clearly, increasing risk implies increasing the variance; however, as described in the introduction, the converse does not always hold: it does for the Gaussian framework, where by construction risk is equivalent to variance, but in all other cases it need not be. Conditions \eqref{MPSCONDO} and \eqref{MPSCONDO2} allow one to generalize the original parameterization to scalar diffusion processes, and thus to a \emph{dynamic} stochastic second-order partial ordering. The first condition (antisymmetry) guarantees the mean-preserving property, in order to keep the ordered processes equivalent in terms of first-order stochastic dominance. The positivity condition is what determines the increase in risk, since for an increase in $\lambda$ the probability weight determined by $\mathcal{P}^{(\lambda)}(x,t)$  at a given point $y < \infty$ increases as well: this is equivalent to the original Rothschild-Stiglitz definition of thicker tails for a riskier process. Finally, note that if $\mathcal{P}^{(\lambda)}(x,t) $ is stopped at an arbitrary time $s\in [0,T]$ then conditions \eqref{MPSCONDO} and \eqref{MPSCONDO2} reduce exactly to \eqref{MPSRS1} and \eqref{MPSRS2}. \\

The following Proposition gives the sufficient condition for a $\lambda$-parameterized distribution to satisfy the conditions \eqref{MPSCONDO} and \eqref{MPSCONDO2}, and therefore be a DMPS:\\

\begin{proposition}\textsc{Sufficient Condition for a DMPS.}\\
Let $R^{(\lambda )}(x):={\frac{\partial }{\partial \lambda }}h^{(\lambda
)}(x)$. \ A sufficient condition for any stochastic process $X_t^{(\lambda)}$ that obeys the TPD \eqref{DUALTPDA} to satisfy the integral conditions (\ref{MPSCONDO}) and (\ref{MPSCONDO2}) is:%
\begin{equation}
R^{(\lambda )}(x)=R^{(\lambda )}(-x)\geq 0.  \label{SCONDO}
\end{equation}
for all $x \in \mathbb{R}$.
\end{proposition}

\vspace{0.3cm}

\begin{proof}
See Appendix \ref{proofprop1}.
\end{proof}
\bigskip
\\
Proposition 1 immediately allows one to characterize dynamic second-order stochastic dominance in terms of the preferences of a risk-averse agent, as stated in Definition 1.4. We do so in the following Lemma: \\
\\
\textbf{Lemma 1: } \emph{For any two ordered stochastic processes $X_t^{(\lambda_1)}$ and $X_t^{(\lambda_2)}$ with $\lambda_1<\lambda_2$ that satisfy the sufficient condition \eqref{SCONDO}, a risk-averse agent with time-consistent and time-invariant preferences will favor $X_t^{(\lambda_1)}$, i.e. $ u(X_t^{(\lambda_1)},t ) \geq  u(X_t^{(\lambda_2)},t )$ for all $t\in [0,T]$,  for all utility functions such that $u_{xx} \leq 0$.  } \\
\\
\begin{proof}
See Appendix \ref{prooflemma1}. 
\end{proof}
\\

Let us now characterize one of the $\lambda$-family diffusion processes that satisfies the antisymmetry and positivity properties of Proposition 1 and therefore is a DMPS. We restrict our attention to processes that exhibit the Brownian bridge property: if a process is conditioned to be 0 at both $t=0$ and $t=T$, then the resulting process is a Brownian bridge. We prove in the following Proposition that for such processes the functional form with non-constant drift that satisfies the conditions in Definition 3 is unique. In the context of economics, such a restriction is not particularly stringent since if one considers a drifted Brownian motion that is conditioned to be 0 at both endpoints, and restricts attention to constant drifts $b =k \in \mathbb{R}$, it is well known that the resulting process is a Brownian bridge. The loss of generality is therefore negligible.\\
%

\begin{proposition}
\textsc{Uniqueness.} \\For any diffusion process which exhibits the Brownian
bridge property, the \textbf{only} diffusion process with non-constant drift which satisfies the DMPS integral conditions \eqref{MPSCONDO} and \eqref{MPSCONDO2} is the diffusion process: 
\begin{equation}\label{prop2}
dX_{t}= \left  [\sqrt{2\lambda}\tanh(\sqrt{2\lambda} X_t )\right  ] dt+dW_{t}.
\end{equation}
with $X\in \mathbb{R}, t\in [0,T], \lambda \in \mathbb{R}^+$.
\end{proposition}
\begin{proof}
See Appendix  \ref{proofprop2}.
\end{proof}
\\ 
%
\\
The process given by \eqref{prop2} is clearly an It\^o process. We now derive its TPD. Its density $Q(x,t|x_0, 0)$ solves the Kolmogorov Forward equation:
\begin{equation}
\frac{\partial}{\partial t} Q(x,t|x_0,0) = -\frac{\partial}{\partial x} \sqrt{2\lambda}\tanh(\sqrt{2\lambda} X_t ) Q(x,t|x_0,0) + \frac{1}{2} \frac{\partial^2}{\partial x^2}Q(x,t,x_0,0).  \label{kolm}
\end{equation}
%
%
%
Solving \eqref{kolm} yields the following Lemma, which will be the workhorse result of the remainder of this paper.\\
\\
\textbf{Lemma 2.} \textsc{Bernoulli Representation Lemma.}\\
\emph{The process given in \eqref{prop2} obeys the following TPD:}

\begin{equation}
Q(x,t|x_0,0) = {\frac{1}{2\sqrt{2\pi t}}}\left\{ e^{-{%
\frac{(x-x_0-\sqrt{2\lambda }t)^{2}}{2t}}}+e^{-{\frac{(x-x_0+\sqrt{2\lambda }t)^{2}}{%
2t}}}\right\} .\label{TPDTANHL}
\end{equation}
\emph{which is a superposition of two $\pm \sqrt{2\lambda }$-drifted Wiener
processes. This implies that equation (\ref{prop2}) can be rewritten as:}%
\begin{equation}
dX_{t}=\mathcal{B}dt+dW_{t},  \label{BERNOULLI}
\end{equation}%
\emph{where $\mathcal{B}$ stands for a Bernoulli random variable taking the values 
$\pm \sqrt{2\lambda }$ with 0.5 probability.} \\
\\
\begin{proof}
See Appendix D.
\end{proof}
\\
\\
Using equation \eqref{TPDTANHL} one can obtain the first moment and the covariance of the process: 
\begin{equation}
\begin{cases}
\mathbb{E} \{ X_t \}  = 0,\\
Cov \{ X_s, X_t \} = 2\lambda  st + \min\{ s,t \} . \label{moments}
\end{cases}
\end{equation}
When $s=t$ one can immediately see that the variance increases quadratically in time. To clarify the representation given in Lemma 2, equation \eqref{BERNOULLI} has to be understood as a process in which at time $t$ one observes the realization of the sign of the drift, and afterwards lets the process behave according to the resulting $\pm \sqrt{2\lambda}$-drifted Brownian motion. This generates never-vanishing correlations in the noise. This representation also corresponds perfectly to Rothschild and Stiglitz's ``addition of noise'' condition for a MPS, where a variable is made riskier by means of the addition of zero conditional mean noise, as in Definition 1.1. \\

Theorem 2 in \cite{machinapratt} allows the construction of such noise in a static setting: it states that if two random variables $\tilde{x}, \tilde{y}$ have the respective cumulative densities $F(.), G(.)$ that satisfy the Rothschild-Stiglitz integral conditions, then one can construct a set of random variables $\tilde{\epsilon}$ with zero conditional (on $\tilde{x}$) mean such that $\tilde{y} =  \tilde{x}+\tilde{\epsilon}$.
Since their theorem is valid for an arbitrary time, it is also applicable in our dynamic framework: setting $\lambda =0$, one has $X_t = W_t$, Gaussian noise with TCD $\mathcal{P}^{(0 )}(x,t)$, which obviously satisfies integral conditions \eqref{MPSCONDO} and \eqref{MPSCONDO2}. 
We have proven in Proposition 2 how a $\lambda$-DMPS process also satisfies the two conditions with the transition cumulative density $\mathcal{P}^{(\lambda )}(x,t)$. Creating a set of $\lambda$-indexed Bernoulli variables $\mathcal{B} = \{ \pm \sqrt{2\lambda} \}$ with probability $0.5$ independent of $X$, with $\lambda \in \mathbb{R}^+$, and calling $\tilde{x}_t = W_t$, we have that $\tilde{y}_t = \mathcal{B} + \tilde{x}_t $ holds for all $t$, because the problem is well-posed and equations \eqref{DUALSDA} and \eqref{DUALTPDA} hold in the entire time domain. 
Because of \eqref{BERNOULLI}, $\tilde{y_t}$ has the TCD given by \eqref{TPDTANHL}. As proven above, $\mathcal{P}^{(0) }$ and $\mathcal{P}^{(\lambda )}$ satisfy the integral conditions for all $\lambda \in \mathbb{R}^+$ and $t \in [0,T]$ and therefore the theorem applies. It can easily be shown that the same applies between any $\lambda$-densities $\mathcal{P}^{(\lambda_1 )}, \mathcal{P}^{(\lambda_2 )}$,  by means of a Bernoulli variable taking values $\pm | \lambda_1 - \lambda_2 | $. \\
\begin{figure}
\centering
\includegraphics[width=14cm]{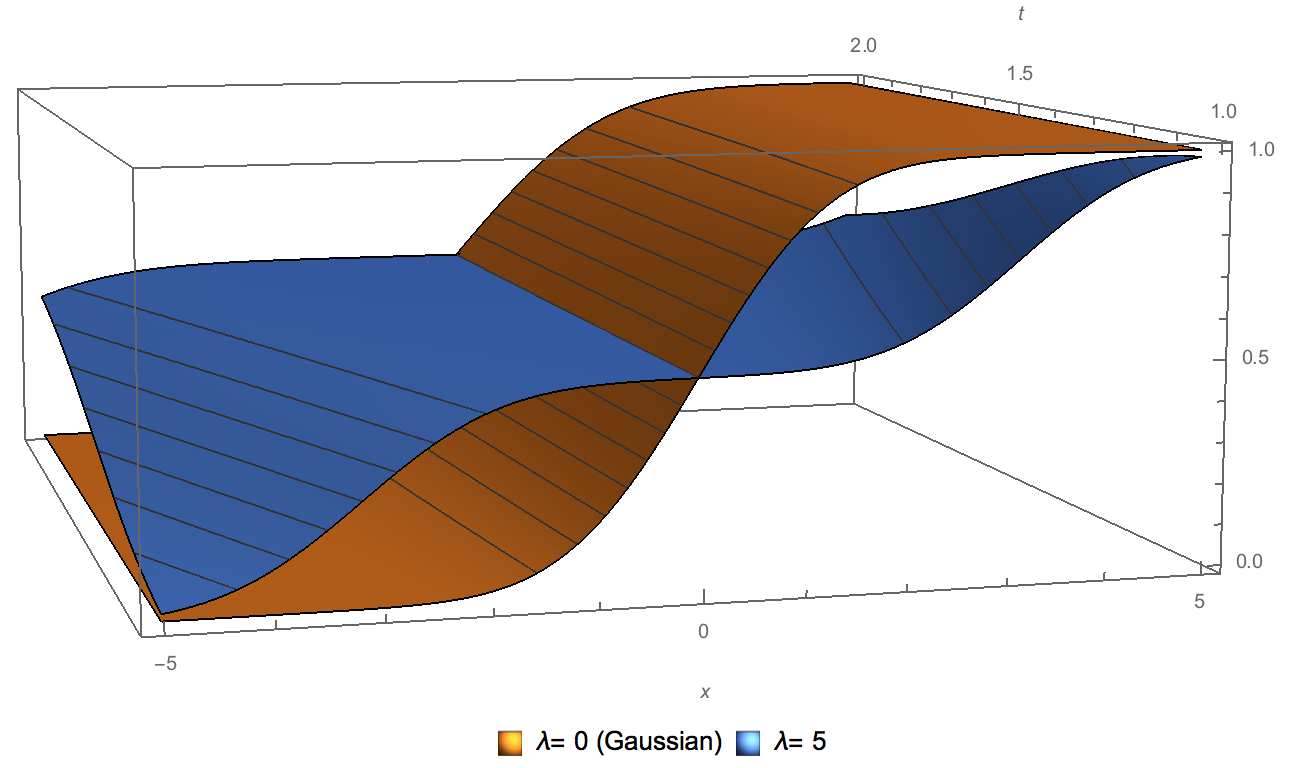}
\caption{An illustration of $\mathcal{P}^{\protect\lambda }(x,t):=\protect\int_{-\infty }^{x}{\frac{1}{2\protect%
\sqrt{2\protect\pi t}}}\left\{ e^{-{\frac{(y-\protect\sqrt{2\protect\lambda }%
t)^{2}}{2t}}}+e^{-{\frac{(y+\protect\sqrt{2\protect\lambda }t)^{2}}{2t}}%
}\right\} dy$, for $\protect\lambda =0$ (Gaussian) and $\lambda=5$. Notice the thicker tails on the left side of the figure.}
\label{figure1}
\end{figure}


The stochastic process given by \eqref{prop2} or equivalently \eqref{BERNOULLI} is of central importance for the economic applications we consider. \cite{Hongler2006} refer to it as \emph{super-diffusive ballistic noise}.  The super-diffusive nature of this process is apparent in \eqref{moments}, in that the variance increases quadratically in time, as previously noted. The two values $\pm \sqrt{2\lambda}$ or, equivalently, the hyperbolic tangent in the drift, is what shifts probability to the tails of the distribution, thereby allowing the process to satisfy the second integral condition \eqref{MPSCONDO2} of the increasing risk definition. \\

\begin{figure}
\centering
\begin{subfigure}[b]{0.8\textwidth}
\includegraphics[width=14cm]{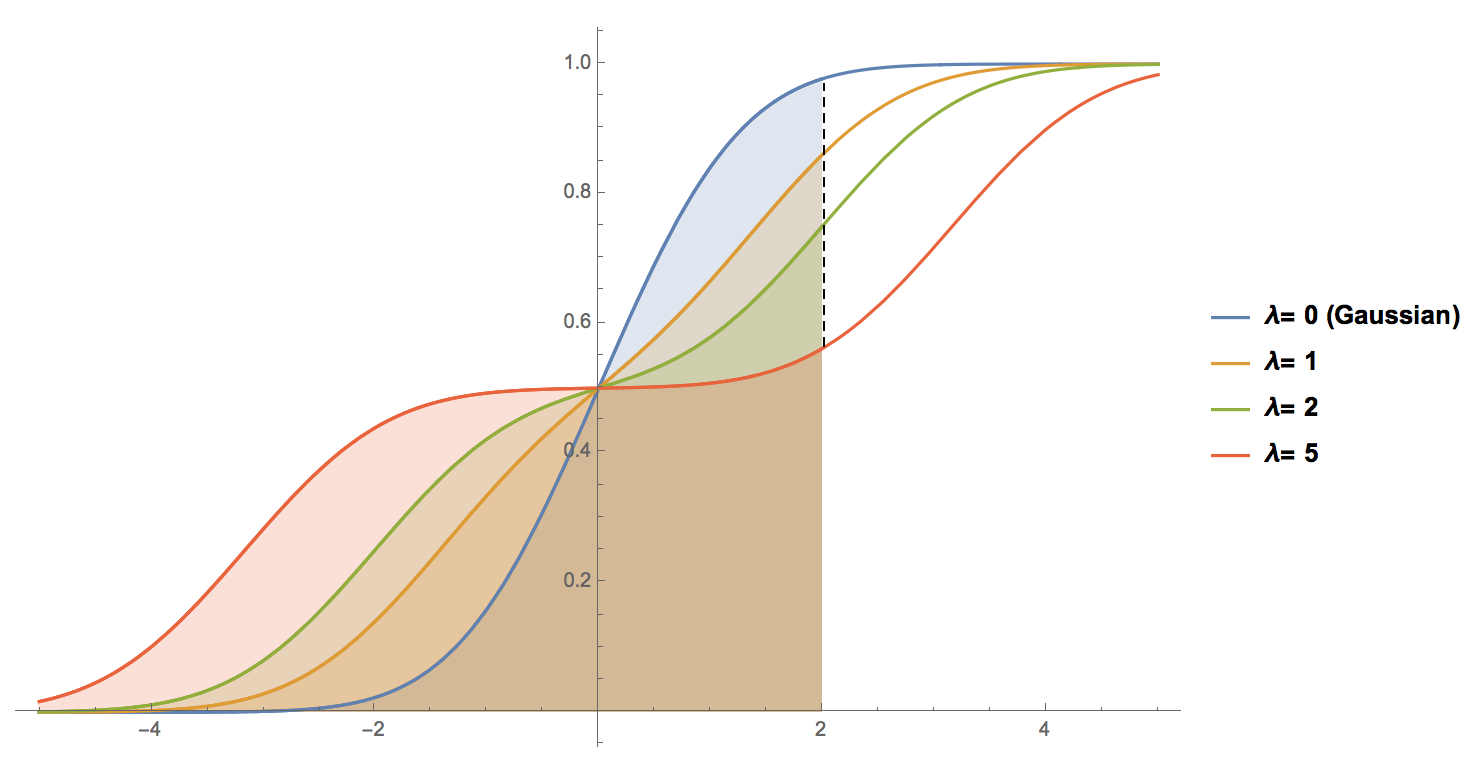}
\caption{An illustration of the positivity condition $(ii)$ in Definition 3 along with changes in concavity and the thicker tails, with $\mathcal{P}^{(\lambda)}(x,t)$ evaluated at $t=1$.}
\label{figure3}
\end{subfigure}\\
\begin{subfigure}[b]{0.8\textwidth}
\includegraphics[width=14 cm]{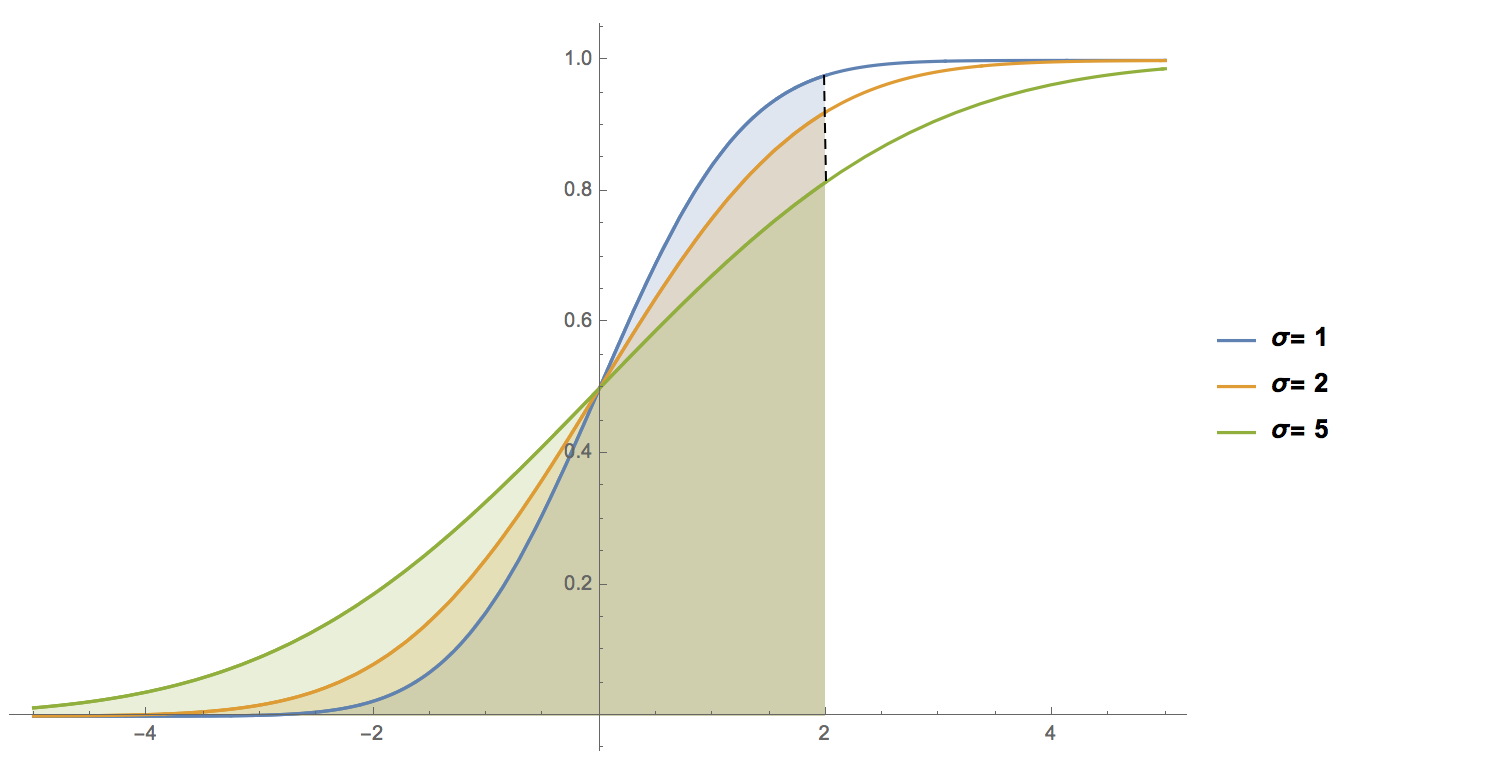}
\caption{An increase in the variance of a driftless Brownian motion $\mathcal{P}^{(0)}(x,t)$ at $t=1$. }
\label{gauss}
\end{subfigure}
\caption{}
\end{figure}

An illustration of $\mathcal{P}^{(\lambda)
}(x,t):=\int_{-\infty }^{x}Q^{(\lambda)}(y,t|x_0,0)dy$, with $Q^{(\lambda)
}(x,t|x_00)$ given by equation (\ref{TPDTANHL}), and for $x_0 = 0, \lambda =0$ (Gaussian) and $\lambda =5$,
is provided in Figure \ref{figure1}. Note the changes in concavity between the Gaussian case and the DMPS. This behavior is structurally impossible to obtain by changing the variance of a Gaussian: in the DMPS process the random superposition of two Gaussian distributions generates thick tails in the distribution while still being able to rank them in order of riskiness. Notice the leftmost portion of Figure \ref{figure1} and observe how the tails thicken as $\lambda$ increases, though the distribution never becomes fat-tailed in a formal sense: the moment generating function $\mathbb{E} e^{\alpha X}$ of the DMPS process is easily shown to be finite for all $0< \alpha < \infty $. This stochastic process has the extremely useful property of escaping the Gaussian framework while remaining tractable analytically. \\

At a more general level, consider two $\mathcal{P}^{(\lambda) }(x,t)$
surfaces for $\lambda_1$ and $\lambda_{2} $, with $\lambda_1 < \lambda_2$. \ For a given $t$, these would
correspond to two cumulative densities in the standard Rothschild and
Stiglitz graphical illustration, with the former having lower risk than the
latter. \ It is clear from Figure \ref{figure3} that the area corresponding to the vertical
distance between the two curves to the right of $x=0$ (the mean) and to the
left of $x=0$ are equal: this corresponds to the first integral condition in
Definition 2. \ If one were to ``stop'' at some value $y$ in $\mathbb{R}^+$, the total ``positive" distance between the two curves over 
$(-\infty $, $0)$ would of   necessity outweigh the ``negative" distance
between the two curves over $(0,y)$: this corresponds to the second integral
condition of Definition 2. This is illustrated in Figure \ref{figure3}, a cross-section of $\mathcal{P}^{(\lambda)}$ at $t=1$ where the dashed line on the right quadrant highlights the positivity. While the point is almost trivial, for direct comparison purposes Figure \ref{gauss} shows the impact of a change in variance in a driftless Brownian motion: one can immediately see that modifying the risk parameter $\lambda$ versus the variance imply very different consequences. Obviously if one sets $\lambda=0$ then one reverts to the pure Gaussian framework, where risk and variance coincide.\\

So as to furnish researchers with the complete panoply of tools allowing them to use the DMPS noise source in economic applications, we now study the behavior of a stochastic process that is driven by equation \eqref{prop2} instead of Gaussian noise. We (i) characterize its probabilistic properties and (ii) derive the appropriate It\^o formula. Define: 

\begin{equation}
\left\{ 
\begin{array}{l}
dZ_{t}=\mu(Z_t,t ) dt+\sigma(Z_t,t ) dX_{t}, \\ 
dX_{t}=\sqrt{2\lambda }\tanh  ( \sqrt{2\lambda }X_{t}) 
dt+dW_{t}, \\
Z_{0}=z_{0}, X_0 = 0.
\end{array}%
\right .
 \label{SCALARMPS}
\end{equation}%
The process \eqref{SCALARMPS} is a degenerate diffusion process in $\{-1, +1\} \times \mathbb{R}^{2}$:
it is characterized by the TPD $P(z,x,t|z_0,0,0)$ 
that solves the Kolmogorov Forward equation: 
\begin{equation}
{\frac{\partial }{\partial t}}P(z,x, t|z_{0},0,0)=\mathcal{F} \left (P(z,x,t|z_{0},0,0)\right ),
\label{FP2}
\end{equation}%
where the operator $\mathcal{F}(.)$ is given by:

\begin{eqnarray}
\mathcal{F}(\cdot ) &=&
-{\frac{\partial }{\partial z}}\left[ \mu(Z_t,t) +\sigma(Z_t,t)  \sqrt{2\lambda } \tanh ( \sqrt{
2\lambda }x ) \right]  -{\frac{\partial }{\partial x}}\left[  \sqrt{2\lambda } \tanh ( \sqrt{2\lambda }x
) \right ] \nonumber \\
&+ &\quad  \overset{=(\partial z,\partial x)\mathbf{\Sigma }\,
\mathbf{\Sigma }^{T}(\partial z,\partial x)^{T}}{\overbrace{\left\{ {\frac{\sigma(Z_t,t)^2
}{2}}{\frac{\partial ^{2}}{\partial z^{2}}}+{\sigma(Z_t,t)\frac{\partial ^{2}}{\partial
z\partial x}}+{\frac{1}{2}}{\frac{\partial ^{2}}{\partial x^{2}}}\right\} }}.
 \label{KOLMOOP}
\end{eqnarray}
In view of equation (\ref{SCALARMPS}), we are interested in the marginal measure of the process $Z_t$, $P_{M}(z,t)$, which results from the $x$-integration: 

$$
P_{M}(z,t)=\int_{\mathbb{R}}P(z,x,t)dx.
$$
For general functional forms $\mu(Z_t,t), \sigma(Z_t, t)$ one cannot solve \eqref{FP2} in closed form. However, we can obtain closed form solutions for three cases of common use in economics: the drifted process with scalar coefficients, the mean-reverting (stationary Ornstein-Uhlenbeck) process and the geometric process. This is done in the following Proposition:

\begin{proposition}\textsc{Marginal densities of DMPS-driven stochastic processes.}\\

3.1 Drifted (scalar) case, $\mu(Z_t,t) = \mu$, $\sigma(Z_t, t)=  \sigma$:

\begin{equation}
P_M(z,t |z_0, 0) = \frac{1}{2 \sqrt{2\pi \sigma^2 t}} \left (e^{-\frac{[z- z_0 -(\mu -\sigma \sqrt{2\lambda})t ]^2}{2\sigma^2 t}} + e^{-\frac{[z- z_0 -(\mu + \sigma \sqrt{2\lambda})t ]^2}{2\sigma^2 t}} \right).\label{prop3}
\end{equation}%

3.2 Mean-reverting case, $\mu(Z_t,t) =   \alpha( \mu -  Z_t)$, $\sigma(Z_t, t) = \sigma$:

\begin{equation}
P_M(z,t |z_0, 0) =  \frac{\sqrt{\alpha}}{2 \sqrt{\pi S(t)}} \left (e^{-\frac{\alpha [ z - z_0e^{-a t} - M^+(t)]^2}{S(t)}} +e^{-\frac{\alpha [ z - z_0e^{-a t} + M^-(t)]^2}{S(t)}} \right), \label{prop3ou}
\end{equation}

where $S(t) =  \sigma^2 (1-e^{-2 \alpha t})$ and $M^{\pm}(t) = (\mu \pm  \sigma \sqrt{2\lambda}/\alpha)(1-e^{-\alpha t})$ .\\

3.3 Geometric case, $ \mu(Z_t,t) =  \mu Z_t$, $\sigma(Z_t, t) = \sigma Z_t $:

\begin{equation}
P_{M}(z,t\mid z_{0},0) = {\frac{1}{2z\sqrt{2\pi \sigma ^{2}t}}}\left( e^{-{
\frac{\left[ \ln (z)-\ln (z_{0})-(\mu -\sigma \sqrt{2\lambda })t\right] ^{2}
}{2\sigma ^{2}t}}}+e^{-{\frac{\left[ \ln (z)-\ln (z_{0})-(\mu +\sigma \sqrt{
2\lambda })t\right] ^{2}}{2\sigma ^{2}t}}}\right), \label{prop3geom}
\end{equation}

for all $\alpha, \mu, \sigma \in \mathbb{R}^+, z \in \mathbb{R}, t \in [0,T]$.
\end{proposition}

\begin{proof}
See Appendix  \ref{pro4}.
\end{proof}
\bigskip
\\
Now consider the simplest example that allows one to jointly use Lemma 1, Proposition 2 and Proposition 3. Consider an individual with CARA utility function of the form $u(x) := -e^{- \gamma x}/ \gamma$, where $ \gamma$ is the Arrow-Pratt coefficient of absolute risk-aversion. Consider an arbitrary process $Z_t$ that evolves as a drifted scalar DMPS process $dZ_t = 
\mu dt + dX_t$, and thus with \eqref{prop3} as its TPD. Let us discretize the dynamics at $t=1$ and assume that $x_0=0$ for simplicity. Because of the special form of the DMPS density we can calculate the agent's expected utility at $t=0$ as a Laplace transform: 
\begin{eqnarray}
\mathbb{E} u(Z_1) &=& -\frac{1}{\gamma} \exp[- \gamma \mu ] \times \nonumber \\
&& \frac{1}{\sqrt{2 \pi}} \int_{-\infty}^\infty  \exp[- \gamma z] \left ( \frac{1}{2} e^{\frac{-(z-\sqrt{2\lambda})^2}{2}} + \frac{1}{2} e^{\frac{-(z+\sqrt{2\lambda})^2}{2}} \right ) dz, \nonumber \\
& =& -\frac{1}{\gamma}  \exp \bigg [- \gamma\left (\mu - \frac{1}{2} \gamma - \frac{1}{ \gamma} \log \cosh ( \gamma \sqrt{2\lambda})\right ) \bigg ] .
\end{eqnarray}
This last expression transparently shows how the expected utility of a risk-averse agent is strictly decreasing in the risk parameter $\lambda$.\\

We now present It\^o's formula to calculate the stochastic differentials of a process driven by a DMPS noise source. Define:

\begin{equation}
\left\{ 
\begin{array}{l}
dZ_{t}=\mu(Z_t,t )  dt+\sigma(Z_t,t ) dX_{t},\qquad Z_{0}=z_{0},\ \, \\ 
dX_{t}=\sqrt{2\lambda }\tanh  ( \sqrt{2\lambda }X_{t} ) 
dt+dW_{t}\qquad X_{0}=0,%
\end{array}%
\right.  \label{SCALARMPS2}
\end{equation}%
where $\mu(.), \sigma(.)$ are real-valued bounded and measurable functions that obey standard conditions for existence and uniqueness of a solution, as given by Theorem 5.2.1 in \cite{Oksendal}, and let $g(Z_t,X_t,t )$ be a real-valued function $g: t \times \mathbb{R}  \to t \times \mathbb{R}$. We then have the following Proposition:

\begin{proposition}\textsc{It\^o's formula for a DMPS process.}\\
Let $g(Z_t, t)$ of class $\mathcal{C}^{2,1}$ on $t \times \mathbb{R} \to t \times \mathbb{R.}$  We then have: 
\begin{eqnarray}
d g(Z_t,t) &=  & \left [ \frac{\partial g}{\partial t} + \bigg ( \mu(Z_t, t)  + \sigma(Z_t, t) \mathcal{B}\bigg) \frac{\partial g}{\partial z} +\frac{\sigma(Z_t,t)^2}{2}\frac{\partial^2 g}{\partial z^2} \right] dt  \nonumber \\
&+&  \sigma(Z_t,t) \frac{\partial g}{\partial z}  dW_t.
\label{itodmps}
\end{eqnarray}%
where $\mathcal{B}$ is defined in the Bernoulli Representation Lemma.
\end{proposition}
\begin{proof}
See Appendix  \ref{itolemma} for the proof and for a more general version of It\^o's formula for a DMPS process when $g(.)$ is also explicitly a function of the noise variable $X_t$, i.e. $g(Z_t, X_t, t)$.
\end{proof}
\newpage 


\section{On the fundamental difference between a DMPS and an increase in variance in a dynamic setting}

We now investigate the consequences of driving a general stochastic process $Z_t$ by means of the DMPS noise source $X_t$ as given by \eqref{prop2} instead of the Brownian motion $W_t$. We do so in order to shed light on the radical difference between a change in variance and a change in risk once one escapes the Gaussian framework. This is done in three parts. 
First, we study the time-invariant (stationary) probability measure of the DMPS-driven process and show how an increase in risk thickens the tails of the distribution and yields bimodality. Second, we study two examples of optimal stopping and show how the continuation decision, as well as the stopping threshold, depends critically on the alternation of the random drift. Third, we show how for a controlled diffusion process driven by DMPS, the Certainty Equivalence Principle is violated and one can no longer use the expectation of the drift to optimally control the process.

\subsection{The curvature of the time-invariant probability measure}

\vspace{0.3cm}

To show that using a DMPS as the driving noise in diffusion processes leads
to behavior that is drastically different from that stemming from an
increase in the variance term $\sigma $ in front of the White Gaussian Noise
(WGN), consider the behavior of the time-invariant (or stationary) probability measure for
scalar processes. Consider the process:

$$
dZ_{t}=f(Z_t)dt+\sigma dW_{t}.  \label{SCALARDIFF}
$$
with $f: \mathbb{R} \to \mathbb{R}$ admitting an antiderivative. This stochastic process is a diffusion process with a general non-constant drift, driven by WGN. The Kolmogorov forward equation for the TPD $P(z,t|z_{0},0)$
associated with equation (\ref{SCALARDIFF}) is:%
$$
{\frac{\partial }{\partial t}}P(z,t|z_{0},0)=-{\frac{\partial }{\partial z}}%
\left[ f(z)P(z,t|z_{0},0)\right] +{\frac{\sigma ^{2}}{2}}{\frac{\partial ^{2}%
}{\partial z^{2}}}P(z,t|z_{0},0),  \label{FOKKER}
$$
and the time-invariant (or stationary) measure $P_{s}(z)=\lim_{t\rightarrow
\infty }P(z,t|z_{0},0)$ is obtained by solving the KFE with a vanishing left-hand side since we have $\partial P_s / \partial t = 0$. Integrating with respect to $z$
(with vanishing constants of integration, since in a stationary state no probability current is sustained), we obtain: %
$$
0=-\left[ f(z)P_{s}(z)\right] +{\frac{\sigma ^{2}}{2}}{\frac{\partial }{%
\partial z}}P_{s}(z).
$$
Integrating again yields:

\begin{equation}
P_{s}(z)=\mathcal{N}e^{{\frac{2}{%
\sigma ^{2}}}F(z)},  \label{PS}
\end{equation}
where $F(x)$ is the antiderivative of $f(x)$ and $\mathcal{N}$ is a normalization
factor which exists for globally attracting drifts (i.e. $%
\lim_{|x|\rightarrow \infty} F(x) = -\infty $). From equation (\ref{PS}), it
is clear that in a Gaussian setting increasing the variance $\sigma ^{2}$ spreads the $P_{s}(x)$
without affecting its extrema.\\

Let us now consider the impact of driving a stochastic process with the DMPS noise source, as done previously in Proposition 3. We have:%
\begin{equation}
\left\{ 
\begin{array}{l}
dZ_{t}=f(Z_t) dt+\sigma dX_{t}, \\ 
dX_{t}=\sqrt{2\lambda }\tanh ( \sqrt{2\lambda }X_{t} )
dt+dW_{t}, \\
Z_{0}=z_{0}, X_0 = 0.
\end{array}%
\right.  \label{SCALARMPS3}
\end{equation}%
The process \eqref{SCALARMPS3} is characterized by the TPD $P(z,x,t|z_0,0,0)$ that solves the Kolmogorov forward equation: 
$$
{\frac{\partial }{\partial t}}P(z,x, t|z_{0},0,0)=\mathcal{F} \left (P(z,x,t|z_{0},0,0)\right ).
$$
where the operator $\mathcal{F}(.)$ is given by \eqref{KOLMOOP}, and we have replaced $\mu$ by $f(z)$.
%
%
%
\begin{figure}
\centering
\includegraphics[width=17cm,scale=1.5]{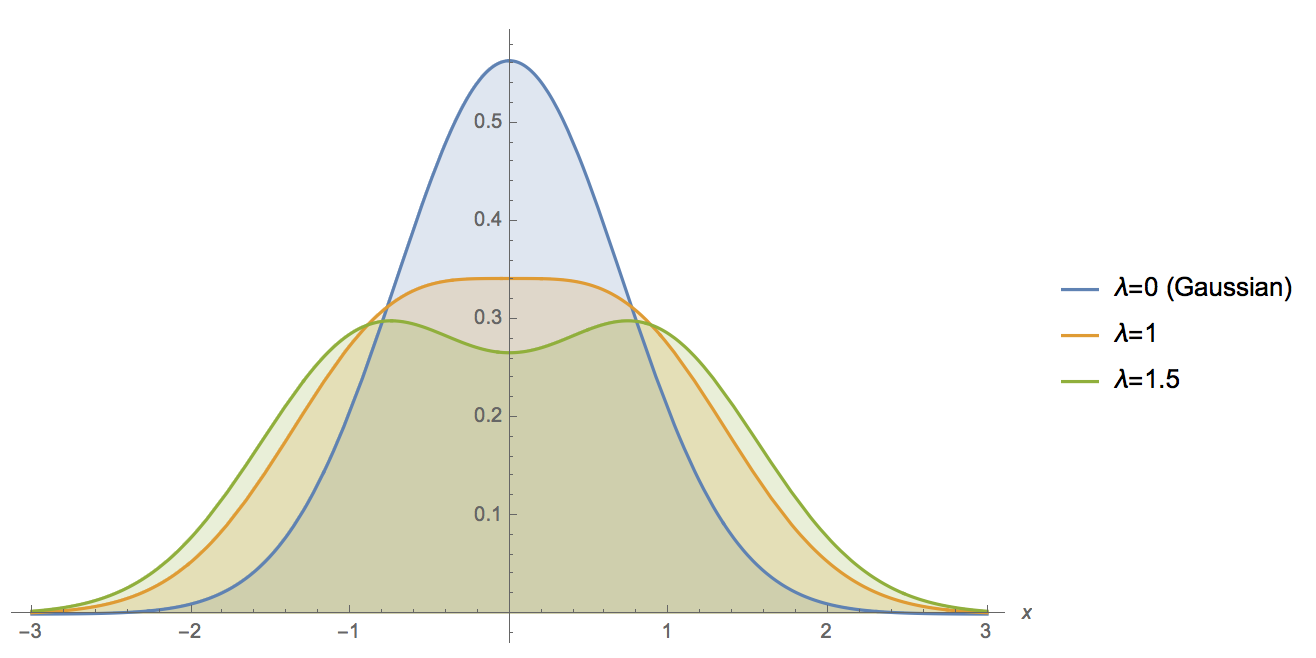}
\caption{Modulating the risk parameter $\lambda$ in the DMPS stationary probability measure. Notice the change in curvature at the origin for $\lambda = 1.5$ and the bimodality of the distribution}
\label{figure2}
\end{figure}
The stationary measure $P_{s}(z,x)$ in $\mathbb{R}^{2}$ then solves $\mathcal{F}%
\left( P_{s}(z,x)\right) =0$. For arbitrary $f(z)$, $P_{s}(z,x)$ and $\mathcal{F%
}\left( P_{s}(z,x)\right) =0$ cannot generally be integrated in closed form.  As done previously for Proposition 3, we are mainly interested in the marginal stationary measure of the DMPS-driven process $P_{s}(z)$, resulting from the $x$-integration of $P_{s}(z,x)$. Application of equation \eqref{BERNOULLI} of the Bernoulli Representation Lemma allows one to rewrite \eqref{SCALARDIFF} as:

$$
dZ_{t}=\left[ f(Z_t)+\mathcal{B}\right] dt+\sigma dW_{t},\qquad Z_{0}=z_{0}.
\label{BERNOUFEED}
$$
Proceeding as in the Gaussian case, one obtains the marginal stationary measure:

\begin{equation}
P_{sM}(z)=\mathcal{N}\left[ e^{\frac{2F^{+}(z)}{\sigma ^{2}}}+e^{\frac{%
2F^{-}(z)}{\sigma ^{2}}}\right] =\mathcal{N}\cosh (\sqrt{2\lambda }z)e^{%
\frac{2F(z)}{\sigma ^{2}}},  \label{PSMB}
\end{equation}%
where $F^{\pm }(z)=F(z)\pm \sqrt{2\lambda }z$. Comparing
equation (\ref{PS}) and equation (\ref{PSMB}), it is obvious that driving
the system with the DMPS process is completely different from simply
modifying the variance. This difference becomes obvious when studying the curvature of the DMPS stationary measure at the origin $\rho
_{(\sigma ,\lambda )}=[ {\frac{\partial ^{2}}{\partial z^{2}}}P_{sM}(z)%
] \mid _{z=0}$, which can be decomposed into two parts. The first corresponds to the Gaussian curvature, while the second depends on the Rothschild and Stiglitz parameter of increasing risk $\lambda$. This is done in the following Proposition:

\begin{proposition}\textsc{Curvature of the stationary DMPS measure.}\\
The curvature at the origin of the stationary marginal measure $\rho
_{(\sigma ,\lambda )}$ of a DMPS-driven stochastic process is given by:
\begin{equation}
\rho _{(\sigma ,\lambda )}(0) = \rho _{(\sigma ,0)}(0) +2\lambda,  \label{CURVA}
\end{equation}
were $\rho _{(\sigma ,0)}$ is the curvature at the origin of the Gaussian probability measure.
\label{prop5}
\end{proposition}

\begin{proof}
See \ref{proofprop5}.
\end{proof}
\\
\\
Proposition 5 shows that the behavior obtained by driving the system with the DMPS process cannot be derived from a simple change of the variance of a Gaussian. The risk parameter $\lambda$ and the variance term $\sigma$ induce very different effects on the stationary probability measure. For example, since $%
\rho _{(\sigma ,0)}(0) <0$ but it may be the case that $\rho _{(\sigma ,\lambda )}(0) >0$, a DMPS may induce a change in the number of modes of $P_{sM}(z)$. We can immediately see that when $\lambda = 0$ one returns to the Gaussian framework, since the probability spread to the tails is zero. Figure \ref{figure2} illustrates the effect of increasing the risk parameter $\lambda$. As $\lambda$ increases from 1 to 1.5, the tails thicken and bimodality emerges.

\subsection{Optimal stopping}
\noindent We now present two simple examples of optimal stopping when a stochastic process is driven by the DMPS noise source. 
\subsubsection{Stopping at the ultimate maximum} 
Consider the following dynamics:

$$
\label{BS}
\left\{
\begin{array}{l}
dY_t = r Y_t dt  \\ 
dX_t = X_t \left[ \mu dt  + \sigma dW_t\right] \\
Y_0 = 0, X_0 = x_0
\end{array}
\right.
$$

\noindent with $r, \mu >0$, where $Y_t$ is a deterministic process and $X_t$ follows a $\mu$-drifted geometric Brownian motion. The initial noise level is normalized to zero for simplicity. For $t\in [0,T]$, define:

$$
\label{RATIO}
P_t :=  {X_t \over Y_t}\qquad {\rm and } \qquad  M_{T}:=\displaystyle  \sup_{t\in [0,T]} P_t.
$$

\noindent Consider a stopping time $\tau \leq T$ with respect to the filtration ${\cal F}_{t\leq T} :=  \sigma \left( X_u, u \in [0,T] \right)$. For $x>0$, we introduce a utility function $U(x)$, by which one would like to find the stopping time $\tau^{*}$ such that:\footnote{This problem is introduced and discussed in a Gaussian setting by \cite{dutoit}  and \cite{shiryaev}.}

\begin{equation}
\label{ULTIMATE}
\mathbb{E}\left\{ \, U  \left ( {P_{\tau^{*}}\over M_T} \right ) \right\} := \displaystyle \sup_{\tau} \mathbb{E}\left \{ \, U \left ( {P_{\tau}\over M_T} \right ) \right\}.
\end{equation}

\noindent In other words, equation (\ref{ULTIMATE}) recommends stopping at the maximal value in the time interval $[0,T]$ (the ``ultimate maximum'').  In the log-utility case $U(x):= \log(x)$, the problem is tractable.\footnote{The problem is solvable, albeit difficult, for other utility functions, but we choose to limit the presentation to the simplest case to highlight the DMPS dynamics. } Indeed, in this case the stopping decision does not depend on $M_T$ since one immediately solves the SDE with $P_t = p_0 e^{ (\mu - r - \frac{1}{2}\sigma^{2})t + \sigma dW_t}$ with $p_0 = x_0/y_0$ and therefore can write:

\begin{eqnarray}
\label{LOGU}
\displaystyle \sup_{\tau}\mathbb{E}\left \{ \, U \left ( {P_{\tau}\over M_T} \right ) \right\} &=& \displaystyle \sup_{t\in [0,T]}\left (\mu - r -  \frac{1}{2}\sigma^{2}\right )t  \nonumber \\
&=& \left\{ \begin{array}{l}  \left (\mu - r - \frac{1}{2} \sigma^{2} \right )T,  \ \        {\rm if } \quad (\mu - r -  \frac{1}{2}\sigma^{2}) >0,  \\ \\ 0 , \qquad {\rm otherwise.}\end{array}\right. 
\end{eqnarray}

\noindent In view of equation (\ref{LOGU}), the optimal decision is either to stop immediately at $t=0$ or to wait until the endpoint $T$. 
\\

Let us now examine the situation when $X_t$ is driven by the DMPS noise source with amplitude $\gamma = \sqrt{2 \lambda}$. Here, the counterpart to equation (\ref{LOGU}) reads:

\begin{equation}
\label{LOGU1}
\displaystyle \sup_{\tau} \mathbb{E} \, U \left\{ {P_{\tau}\over M_T}\right\}= \displaystyle \sup_{t\in [0,T]} \underbrace{ \left (\mu + {\cal B} \gamma - r -  \frac{1}{2}\sigma^{2} \right )}_{:= \nu_{{\cal B}}}t
=   \displaystyle \sup_{t\in [0,T]} \nu_{{\cal B}} t, 
\end{equation}

\noindent where ${\cal B} = \left\{ -1, +1\right\}$ is the Bernoulli random variable.  In this case we see that the optimal decision is more involved, since besides the two choices in equation (\ref{LOGU}), we have an additional range of possibilities given that the sign of the drift $\nu_{{\cal B}}$ can alternate. In this situation, taking the optimal decision requires additional information, namely the sign of ${\cal B}$ at time $t=0$: once this is obtained, as shown in the previous section, the process becomes Markovian. Again, this shows clearly that driving a stochastic variable with the DMPS leads to behavior which is different from simply increasing the variance of a Brownian motion.

\subsubsection{Stopping with a transactions cost}
Let us now consider another simple illustration, extending Example 10.2.2 by \cite{Oksendal} to our framework. Consider the following dynamics:

$$
\begin{cases}
dX_t =  ( \mu + \gamma \mathcal{B} )X_t dt + \sigma X_t  dW_t , \\
 X_0 = x_0.
\end{cases}
$$
Let us consider the simplest time-inhomogeneous case of optimal stopping, where $c$ is a cost parameter:

$$
f(t,X_t)= e^{-\rho t} (X_t - c).
$$
The stopping problem is to find $\tau^*$ s.t.  $\max_\tau \mathbb{E}_t f(X_\tau)$, for $\tau \geq t$. The value function is given by:

$$
V(t, X_t)  = \sup_{\tau} \mathbb{E}_t \big \{ f(\tau, X_\tau ) | X_t \big \} .
$$
For $g : \mathbb{R}^2 \to \mathbb{R}$, the generator of $X_t$ is:

$$
\hat{\mathcal{A}} g(t, X_t) =  \frac{\partial g}{\partial s} + X_s (\mu + \mathcal{B}\gamma) \frac{\partial g }{ \partial x} + \frac{1}{2} \sigma^2 X_s^2 \frac{\partial^2 g}{\partial x^2},
$$
and the continuation region for the value function therefore becomes

$$
A := \{ (t,x); \hat{\mathcal{A}} V(t,x) >0 \} =
  \begin{cases}
  \mathbb{R} \times \mathbb{R}^+ & \text{if  } \mu + \gamma \mathcal{B} \geq \rho, \\
  x < \frac{c \rho}{\rho - \mu - \gamma \mathcal{B}} & \text{if  } \mu + \gamma \mathcal{B} < \rho.
  \end{cases}
$$
The continuation region does not depend on $\sigma$, but does depend on the risk parameter $\gamma$. The region $A$ is random, since it depends upon the realization of the Bernoulli variable. If $\mu$ is large enough, i.e. $\mu - \gamma >\rho$, then $\tau^* = \infty$. One never stops and a finite optimal stopping time does not exist. Let us assume from now on that this is not the case.  We can then rewrite things as the following boundary value problem:
\begin{eqnarray*}
&&0 = \frac{\partial V}{\partial s} + X_s (\mu + \mathcal{B}\gamma) \frac{\partial V }{ \partial x} + \frac{1}{2} \sigma^2 X_s^2 \frac{\partial^2 V}{\partial x^2},  \\
&& V(s,x_0) = e^{-\rho s}(x_0 - c).
\end{eqnarray*}
Note again that the process $X_t$ itself is not Markovian, because of the correlations generated by $\mathcal{B}$. As such one needs to gather information on the realization of the Bernoulli variable before an optimal decision can be made. 
%
%
By standard smooth fitting arguments, we obtain the optimal stopping rule: 
\begin{equation}
\tau^* \quad  \text{s.t.} \quad X_{\tau^*} = c \frac{\alpha_1(\mathcal{B}) }{\alpha_1(\mathcal{B})  - 1}, \label{stop} 
\end{equation}
where $\alpha_1(\mathcal{B})$ is given by:

$$
\alpha_1(\mathcal{B})  = \frac{\frac{\sigma^2}{2} - \mu -  \gamma\mathcal{B} + \sqrt{\left (\mu + \gamma\mathcal{B} - \frac{\sigma^2}{2}\right)^2 + 2 \rho \sigma^2 }}{\sigma^2}.
$$
We see that the stopping threshold is now a random variable. We also see how the variance of the Brownian motion $\sigma^2$ affects only the stopping threshold \eqref{stop}, while the risk parameter $\gamma$ impacts both the stopping region and the threshold. 

\subsection{A DMPS may violate the Certainty Equivalence Principle}
\vspace{0.3cm}
As in section 2.6 of \cite{Karatzas1997a} consider, for time $t\in
\lbrack 0,T]$, the controlled scalar stochastic process in finite time
horizon defined by the diffusion process:

$$
dZ_{t}=\pi (t)\left[ bdt+\sigma dW_{t}\right] ,\qquad Z_{0}=z\in \lbrack
0,1],  \label{KARA1}
$$
where $b\in \mathbb{R}$ is a constant drift, $dW_{t}$ is the standard Wiener Process and $%
\pi (t)$ is a control. We consider the class $\mathcal{H}(x)$ of admissible
controls $\pi (t)$ which are progressively measurable and for which we have:%

$$
\left \{
\begin{aligned}
&\int_{0}^{T}\pi ^{2}(s)ds<\infty,  \\ 
& 0\leq Z_{t}\leq
1 \ \ \   \mathrm{for} \ \ \in \lbrack 0,T].  \label{ADMI}
\end{aligned}
\right . 
$$
Theorem 2.6.4 of \cite{Karatzas1997a} computes the explicit value
function $G(x)$ of the problem:

\begin{equation}
G(x):=\sup_{\pi (\cdot )\in \mathcal{H}(x)}\mathrm{Prob}\left\{ Z^{x,\pi
}(T)=1\right\} ,  \label{GG}
\end{equation}%
and provides the optimal process $\hat{\pi}(t)$ that attains the supremum in
equation (\ref{GG}). In words, the goal is to determine the optimal control $%
\hat{\pi}(t)$ that maximizes, over the finite time horizon $T$, the
probability of reaching the right-hand boundary $1$ without touching the
left-hand boundary $0$. \\

Now consider the same problem when the noise source $dW_{t}$ is replaced by
the DMPS process of Proposition 2. The problem now requires additional information on the realization of the Bernoulli variable $\mathcal{B}$ in the DMPS process: intuitively, because of the fact that now the noise source can add to the drift $b$ an extra element that may be either positive or negative ($\pm \sqrt{2\lambda}$), one has to consider the possibility of the overall deterministic part of $dZ$ being negative, which was not possible in the original Gaussian formulation. We can rewrite the dynamics as:%
$$
\left\{ 
\begin{array}{l}
dZ_{t}=\pi (t)\left[ bdt+dX_{t}\right] ,\qquad Z_{0}=z\in \lbrack 0,1], \\ 
dX_{t}=\left[ \sqrt{2\lambda }\tanh ( \sqrt{2\lambda }X_{t} ) %
\right] dt+dW_{t}.%
\end{array}%
\right.  \label{KARA2}
$$
Again using (\ref{BERNOULLI}), by the Bernoulli Representation Lemma, we can rewrite the dynamics of $Z_t$ as a random-drifted process:%
\begin{equation}
dZ_{t}=\pi (t)\left[ \hat{b}dt+dW_{t}\right] ,\qquad Z_{0}=z\in \lbrack 0,1],
\label{KARA3}
\end{equation}%
where $\hat{b}$ is a Bernoulli random variable with $\mathrm{Prob}\left\{ 
\hat{b}=b\pm \sqrt{2\lambda }\right\} =1/2$ and hence $\hat{b}$ is drawn from the probability
density function $p(x)dx$:%
\begin{equation}
p_{b}(x)dx={\frac{1}{2}}\left[ \delta (x-b-\sqrt{2\lambda })+\delta (x-b+%
\sqrt{2\lambda })\right] dx,  \label{HATB}
\end{equation}%
where $\delta (x-z)dx$ is the Dirac mass at $z$. Using a martingale
approach, \cite{Karatzas1997} establishes that, provided the support
of $p_{b}(x)$ lies strictly in $\mathbb{R}^{+}$ or in $\mathbb{R}^{-}$, the
optimal control in the presence of a random drift $\hat{b}$ can be directly
obtained from the deterministic case by a simple substitution $b\mapsto \hat{%
b}(t)=\mathbb{E} \{ \hat{b}\mid \mathcal{F}(t) \} $, where $\hat{b}%
(t)$ is the conditional expectation of $\hat{b}$ given the observation of
the process up to time $t$. This is called the Certainty Equivalence
Principle (CEP). Conversely, for cases where the support $p_{b}(x)$ crosses
the origin, the CEP is violated and the resulting optimal control is also
explicitly calculated in \cite{Karatzas1997}. Clearly, in the Gaussian case, one can write $%
p_{b}(x)=\delta (b-x)$ and, for $b\neq 0,$ the support of $p_{b}(x)$ never
crosses the origin. Then for all values of $\sigma $, the CEP holds. This is
not the case for the DMPS process of equation (\ref{KARA3}). Here, when $\sqrt{2\lambda }>b$, the support of $p_{b}(x)$ is
simultaneously contained in $\mathbb{R}^{+}$ and  $\mathbb{R}^{-}$: hence, the CEP does not hold and the optimal control cannot be obtained by a simple substitution of $\hat{b}$. This clearly shows that a DMPS is not equivalent to a modification of the variance parameter $\sigma $
lying in front of a WGN.


\section{Illustrations}

We now present a reworking of four classical economic problems in which the noise that drives the system is our DMPS process instead of WGN: portfolio choice, investment under uncertainty, asset dynamics and entry and exit decisions under uncertainty. Our goal in this section is not to provide novel theoretical models, although some interesting new insights do emerge. Rather, our aim, within the context of four standard economic models, is to highlight both the applicability and the tractability of the DMPS as a tool for applied modelling, while correctly parameterizing the riskiness of the underlying distribution, and clearly distinguishing between risk and variance. 

\subsection{Portfolio selection}

Consider the simplest problem of portfolio selection, in which a risk-averse agent allocates her wealth between a risky asset and a riskless asset with zero return.\footnote{We assume zero return for simplicity, as in a non-interest bearing bank deposit: adding a nonzero interest rate changes nothing.} For simplicity of exposition we will assume CARA preferences, although what follows can be equally solved with CRRA and power utilities. We assume that the return to the risky asset follows a geometric DMPS process given by:

$$
\left\{ 
\begin{array}{l}
dS_t = \mu S_t dt + \sigma S_t dX _{t}, \\ 
dX _{t}=\sqrt{2\lambda} \tanh (\sqrt{2\lambda} X_t)dt+dW_{t}.%
\end{array}%
\right. 
$$
The agent allocates a share $0 \leq v^*(t) \leq 1$ of her wealth to $S_t$ in order to maximize her terminal utility of wealth at a fixed time $T$, which we denote by $u(a_T)$. We model $u(.)$ as a CARA utility function that reads:

$$
u(a_t) = - \frac{1}{\gamma} \exp (-\gamma a_t) ,
$$
where $\gamma$ is the constant Arrow-Pratt coefficient of risk-aversion. If we allow no borrowing ($a_t \geq 0$ for all $t\in [0,T]$), the agent's problem is to find the investment strategy that maximizes her expected utility at time $T$ constrained by a controlled diffusion process.  The optimization can then be written as:

$$
\left\{ 
\begin{array}{l}
\max \mathbb{E}\big [ u(a_T) \big ], \\
\text{s.t. } da(t) = a(t) v(t) \mu dt + a(t) v(t) \sigma dX_t, \\
\quad  \quad  a(0) = a_0 \geq 0.
\end{array}%
\right. 
$$
This is a standard stochastic control problem.  We follow \cite{Oksendal} who, among many others, provides the usual treatment for the Gaussian case. Define a performance function of the form:

$$
J(v(.); t, a) = \mathbb{E}_a \big [u(a_T) \big ]
$$
where $\mathbb{E}_a [ . ]$ denotes conditioning the expectation on $a(t) = a$. Then the problem reduces to finding a Markov control $v^*(t) = v^*(t, a(t))$ such that the individual maximizes the value function:

$$
V(t, a) := \sup_{v(.)} J(v(.); t, a) .
$$
Using the Bernoulli Representation Lemma and Proposition 4, the Hamilton-Jacobi-Bellman (HJB) equation for this problem reads:

$$
\sup_{v} \left \{ \frac{\partial V}{\partial t} +  a v( \mu + \sqrt{2\lambda}\sigma \mathcal{B} )  \frac{\partial V}{\partial a} + \frac{1}{2} a^2 v^2 \sigma^2  \frac{\partial^2 V}{\partial a^2} \right \} =0,
$$
where $\mathcal{B}$ is a $\pm 1$ Bernoulli variable with probability 0.5, and the optimal control is: 

$$
v^*(t,a) =  \frac{ \mu + \sqrt{2\lambda}\sigma \mathcal{B} }{\sigma^2 a} \left(-\frac{V_a}{V_{aa}}\right).
$$
We can already notice that since only the coupled process $\{a_t, X_t\} $ is Markovian, the optimal control is a random (non-Markovian) control until the Bernoulli variable is realized. Once the $\pm 1$ is observed, the optimal control follows immediately. Substituting $v^*$ into the HJB equation yields the following boundary value problem:

\begin{eqnarray*}
& & V_t  - \frac{1}{2} \frac{\mu + \sqrt{2 \lambda} \sigma \mathcal{B}}{\sigma^2}\frac{V_a^2}{V_{aa}}  =0, \\
& &  V(T, a) = -\frac{1}{\gamma} \exp (-\gamma a). \label{portfo3}
\end{eqnarray*}
Using the boundary condition as the basis for a guess for the value function yields:

$$
V(t,a) = -\frac{1}{\gamma} \exp \left [-\gamma \left ( \frac{(\mu + \sqrt{2 \lambda} \sigma \mathcal{B})^2}{2 \gamma} (T-t) + a \right ) \right ].
$$
Once solved for the value function, the optimal control then reads:

$$
v^*(t,a) =  \frac{ \mu + \sqrt{2\lambda}\sigma \mathcal{B} }{\gamma \sigma^2 a}.
$$
This implies that, in contrast to the Gaussian case, the optimal amount of wealth invested in the risky asset can either increase or decrease depending on the realization of the Bernoulli variable, and can potentially be pushed outside the unit interval, in which case the individual would either invest entirely in the risky asset $(v^*=1)$ or deposit everything in the bank $(v^*=0)$. A straightforward application of Proposition 4 to the square of the optimally controlled wealth $a^{*2}_t$ allows one to show that the expected terminal wealth and its variance read:

\begin{eqnarray}
\mathbb{E} (a_T)& =& a_0 + \frac{(\mu +\sqrt{2 \lambda} \sigma \mathcal{B})^2 T }{
\sigma^2 \gamma}, \\
Var (a_T) &=& \frac{(\mu +\sqrt{2 \lambda} \sigma \mathcal{B})^2 T }{(\sigma \gamma)^2}.
\end{eqnarray}
This shows how an increase in risk increases both expected terminal wealth and its variance.  To see this more clearly, simply set $\mu =0$ and one gets $\mathcal{B}^2 = 1$. If one sets $\lambda = 0$ and reverts to a Gaussian setting, an increase in variance $\sigma^2$ decreases both. Conversely, in a DMPS setting, an increase in $\lambda$ (risk) increases both expected terminal reward and riskiness.\footnote{This problem can be equivalently solved, albeit in a slightly more involved way, for CRRA utility functionals: the conclusions in this case retain the same properties as the CARA case. }

\subsection{Optimal investment under uncertainty}

\vspace{0.3cm} 
Following \cite{Abel1983} and \cite{Abel1994}, we consider the
dynamics given by:%
%

%
$$
\left\{ 
\begin{array}{l}
dK_{t}=\left[ I_{t}-\delta K_{t}\right] dt, \\ 
d\epsilon _{t}= \sigma \epsilon _{t} dZ _{t}, \\ 
dZ _{t}=\sqrt{2\lambda} \tanh (\sqrt{2\lambda} Z_t)dt+dW_{t}.%
\end{array}%
\right.  \label{ABELMOD}
$$
where $I_{t}$ stands for investment, $K_{t}$ is the capital stock, $\delta $ is
the depreciation rate, $\lambda$ is the parameter of increasing risk, $\epsilon _{t}$ stands for a stochastic productivity shock and $\sigma $ is a multiplicative noise amplitude factor. For $\mu = 0$ this boils down to the original \cite{Abel1994} formulation. Consider a risk-neutral firm which chooses, over an infinite time horizon,
to maximize the expected present value of operating profit (including investment costs) $\pi \left(
K_{t},I_{t},\epsilon _{t}\right) $:%
$$
V\left( K_{t},\epsilon _{t}\right) =\sup_{I_{t+s}}\int_{t}^{\infty }\mathbb{E%
}\left\{ \pi \left( K_{t},I_{t},\epsilon _{t}\right) \right\}
e^{-rs}ds,  \label{MAXIMAL}
$$
where $r>0$ stands for the discount rate and $\mathbb{E}_{t} \{ \cdot
\} $ stands for the expectation operator conditioned at time $t$. The
corresponding Bellman equation is then given by:%
$$
rV\left( K_{t},\epsilon _{t}\right) =\max_{I_{t}}\left\{ \pi \left(
K_{t},I_{t},\epsilon _{t}\right) +{\frac{1}{dt}}\mathbb{E}\left\{
dV_{t}\left( K_{t},\epsilon _{t}\right) \right\} \right\} .
\label{BELLMAN}
$$
Using Proposition 4 and the Bernoulli Representation Lemma, the Bellman equation can be written as: 
$$
\left\{ 
\begin{array}{l}
rV\left( K_{t},\epsilon _{t}\right) =\max_{I_{t}}\left\{ \pi \left(
K_{t},I_{t},\epsilon _{t}\right) +q_{t}\left( I_{t}-\delta K_{t}\right) +%
\mathcal{L}_{\epsilon }\left[ V\left( K_{t},\epsilon _{t}\right) \right]
\right\} , \\ 
\mathcal{L}_{\epsilon }\left[ V\left( K_{t},\epsilon _{t}\right) \right] :=%
\left[ (\epsilon_t + \mathcal{B} ) {\frac{\partial }{\partial \epsilon }}+{%
\frac{\sigma ^{2}\epsilon_t }{2}}{\frac{\partial ^{2}}{\partial \epsilon ^{2}%
}}\right] V\left( K_{t},\epsilon _{t}\right) ,%
\end{array}%
\right.  \label{MAXDIFF}
$$%
where $q := V_K$ is the marginal valuation of an unit of installed capital and $\mathcal{B}$ is a Bernoulli variable taking values $\pm  \sqrt{2\lambda}$ with probability 0.5. The first order condition for optimal investment $q - \pi_I(K, I^*,\epsilon) = 0$ holds deterministically. We implement for simplicity a Cobb-Douglas production function for a firm that uses labor $L_t$, pays a fixed wage $\omega \geq 0$ and sells its output at a price $P$. This can be made stochastic relatively easily but with no added value for the DMPS illustration so for this purpose we remain in the original framework. Defining $p_t := P \epsilon_t$ and we normalize $P$ to one for simplicity. Abstracting from investment costs (assuming adjustment costs to be independent of the capital stock), we obtain the profit of the optimizing firm: 
$$
\pi \left( K_{t},p_{t}\right) = \max_{L_{t}}\left\{p_t L_t^{\alpha }K_t^{1-\alpha
}-\omega L_{t}\right\}  = h  p_t ^{\theta }K_{t},
$$
with $h:=(1-\alpha )\,\alpha ^{{\frac{\alpha }{1-\alpha }}}\,\omega ^{-{\frac{%
\alpha }{1-\alpha }}}$ and $\theta = 1/(1-\alpha)$ resulting from instantaneous profit maximization with respect to $L_t$. The present value $q_{t}^{(\lambda)}$ at time $t$ of
marginal profits of currently installed capital (taken from time $0$), for a specific risk parameter $\lambda$, is then given by:%
$$
q_{t}^{(\lambda)}=h\int_{0}^{\infty }\mathbb{E}_t\left\{ p_{t+  s}^{\theta
}\right\} e^{-(r+\delta )s}ds.  \label{ACCR}
$$%
Using the Bernoulli Representation Lemma we can write:

$$
d\epsilon_t = \mathcal{B} \epsilon_t dt + \sigma \epsilon_t dW_t,
$$
and using Proposition 3.3 we obtain:

\begin{eqnarray}
q_{t}^{(\lambda)} &=&{\frac{1}{2}}h\,p_{t}^{\theta }\left\{ {%
\frac{1}{r+\delta -\theta\sigma \sqrt{2\lambda} -{\frac{1}{2}}\theta (\theta -1)\sigma ^{2}}}+{\frac{1%
}{r+\delta +\theta \sigma\sqrt{2\lambda} -{\frac{1}{2}}\theta (\theta -1)\sigma ^{2}}}\right\},
\nonumber \\
&=&h\,p_{t}^{\theta }{\frac{(r+\delta )-{\frac{1}{2}}\theta
(\theta -1)\sigma ^{2}}{\left[ (r+\delta )-{\frac{1}{2}}\theta (\theta
-1)\sigma ^{2}\right] ^{2}-2 \theta^2\sigma^2\lambda}}. \label{FINALBALLISTIC} 
\end{eqnarray}
When $\lambda = 0$ (and thus $d\epsilon _{t}=\sigma \epsilon _{t}dW_{t}$), we obtain:

\begin{equation}
q_{t}^{(0)}={\frac{hp_{t}^{\theta }}{%
r+\delta  -{\frac{1}{2}}\theta (\theta -1)\sigma ^{2}}},
\label{FINALMOD}
\end{equation}
which corresponds exactly to the original \cite{Abel1994} result in their equation (30), p. 1379.\\ 

Note that (\ref{FINALBALLISTIC}) cannot be obtained from equation (\ref%
{FINALMOD}) by simply increasing the variance. This illustrates the fact
that a DMPS of a TPD is not equivalent to an increase in the
variance of the normal density. For $\lambda$ relatively small (note that $\theta <1$), equation (\ref{FINALBALLISTIC})
can be approximately rewritten as:

\begin{eqnarray}
q_{t}^{(\lambda)} &\simeq& {\frac{h\,p_{t}^{\theta }}{(r+\delta )-{
\frac{1}{2}}\theta (\theta -1)\sigma ^{2}}}\left[ 1+{\frac{2 \theta^2\sigma^2  \lambda}{\left[
(r+\delta )-{\frac{1}{2}}\theta (\theta -1)\sigma ^{2}\right] }}\right] +
\mathcal{O}(\lambda^2),\nonumber \\ 
\quad &=&q_{t}^{(0)}\left[ 1+2 \lambda \left[ \frac{\theta\sigma  \,q_{t}^{0}}{h\,p_{t}^{\theta }}
\right] ^{2}\right] +\mathcal{O}(\lambda^2).  \label{FINALBALLISTICAPP}
\end{eqnarray}
From equations (\ref{FINALBALLISTIC}) and (\ref{FINALBALLISTICAPP}), we
conclude that, for a given $p_{t}$, the presence of ballistic noise leads to
an increase in the marginal profits $q_{t}^{(\lambda)}$ yielded by
currently installed capital. This is coherent with the original Abel and Eberly framework in which riskiness increases $q$, but cannot be obtained by simply modifying the variance $\sigma^2$ of the Gaussian. In their framework this result is clearly a consequence of the assumption of risk-neutral firms, since an increase in risk means an increase in potential rewards, and riskiness and variance coincide because of their Gaussian setup. 
%

\subsection{Asset dynamics}

\vspace{0.3cm} 

Consider the Black-Scholes asset dynamics $S_t\geq 0$ driven by the
DMPS source $Z_{t}$, namely:%
$$
\left\{ 
\begin{array}{l}
dS_t = \mu S_t \,dt+\sigma S_t dZ_{t}, \\ 
dZ_{t}=\sqrt{2\lambda }\tanh ( \sqrt{2\lambda }Z_{t} ) dt+dW_{t},
\\ 
S_{0}=s_{0}, Z_{0}=z_0,%
\end{array}%
\right.  \label{BSDYN}
$$%
with $\mu, \sigma \in \mathbb{R}^+$. 
%
%
%
%
$$
P_{M}(s,t\mid s_{0},0) = {\frac{1}{2}}P_{\mathrm{BS}}^{(-)}(s,t\mid s_{0},0)+{
\frac{1}{2}}P_{\mathrm{BS}}^{(+)}(s,t\mid s_{0},0),
\label{COSSUM1}
$$
which is the superposition of a pair of Black-Scholes log-normal probability densities $%
P_{\mathrm{BS}}^{(\pm )}(s,t\mid s_{0},0)$ with rates $\mu \pm \sigma \sqrt{%
2\lambda }$. It is worth pointing out that it can be the case that the
ballistic component leads to $\mu -\sigma \sqrt{2\lambda }<0,$ thus
exhibiting a net drifting tendency to be absorbed in the bankrupt state $x=0$ even when the growth rate of the asset is positive, i.e. $\mu >0$ %
. Let us calculate the first moment $\mathbb{E}^{(\lambda )}(S_t\mid
s_{0}) $ for the asset dynamics driven by a DMPS noise source: 
%
%
\begin{figure}
\centering
\includegraphics[width=12cm]{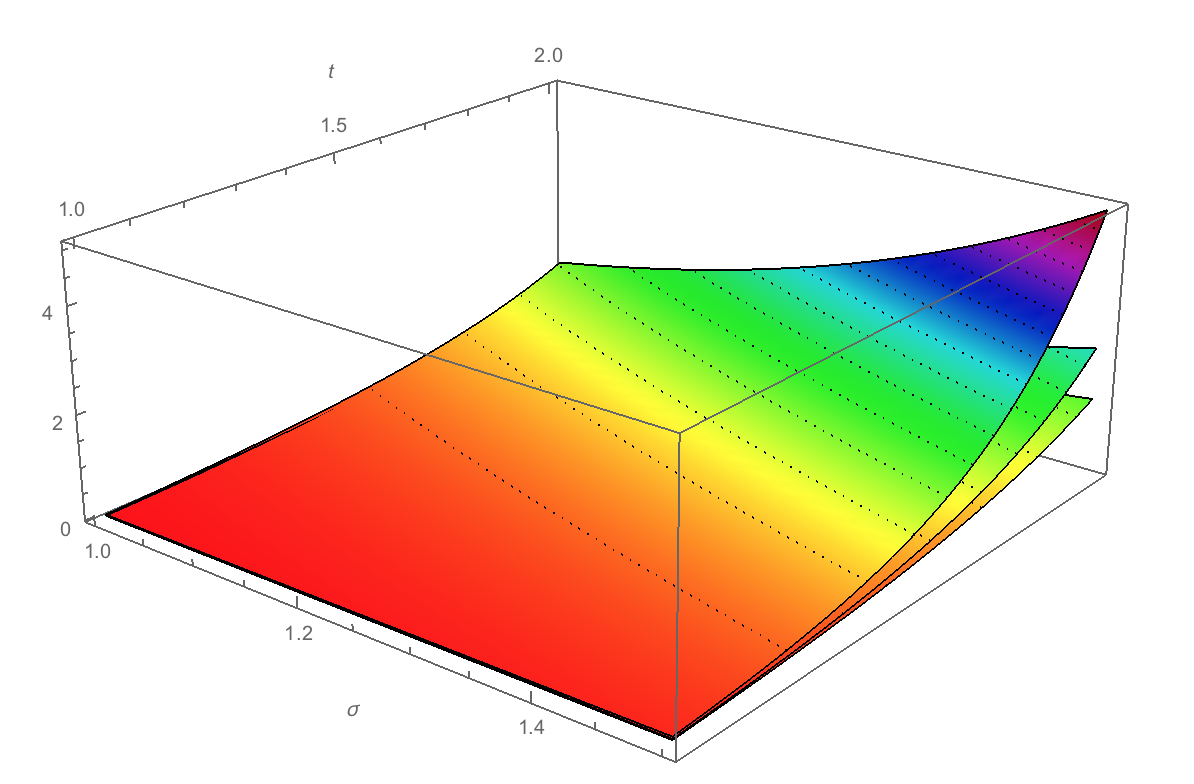}
\caption{An illustration of the
Black-Scholes average asset dynamics driven by DMPS $\mathbb{E}^{(\protect\lambda %
)}(S_t\mid s_{0})=s_{0}e^{\left( \protect\mu +\protect\sigma ^{2}\right)
t}\cosh \left[ \protect\sigma \protect\sqrt{2\protect\lambda }t\right] $,
for $\protect\lambda =0,0.1,0.3$ (we set $s_{0}=0.01,\protect\mu =1$).}
\label{bs1}
\end{figure}
%
\begin{eqnarray}
\mathbb{E}^{(\lambda )}\{ S_t\mid s_{0} \} &=&{\frac{s_{0}}{2}}\left\{
e^{\left( \mu -\sigma \sqrt{2\lambda }+\sigma ^{2}\right) t}+e^{\left(
\mu +\sigma \sqrt{2\lambda }+\sigma ^{2}\right) t}\right\},\nonumber \\
& =& s_{0}e^{\left(
\mu +\sigma ^{2}\right) t}\cosh \left[ \sigma \sqrt{2\lambda }t\right] .
\label{MOMUN}
\end{eqnarray}%
Comparing the average asset growth for WGN with respect to the ballistic
driving noise, one notes that:%
\begin{equation}
{\frac{\mathbb{E}^{(\lambda )}\{S_t\mid s_{0} \} }{\mathbb{E}^{(0)} \{ S_t\mid
s_{0} \} }}=\cosh \left[ \sigma \sqrt{2\lambda }t\right],  \label{COMPARE}
\end{equation}%
and therefore a net average growth enhancement due to the ballistic driving
environment. \ An illustration of equation (\ref{MOMUN}) is provided in
Figure \ref{bs1} for the standard Black-Scholes case ($\lambda =0$), and
progressively higher values of $\lambda $.\\
\\%

\subsection{Entry and exit decisions under uncertainty}

\vspace{0.3cm}

We consider the celebrated \cite{Dixit1989} model, in which a firm undertakes a single discrete project subject to sunk investment costs $k$,
no physical depreciation, and an avoidable operating cost $w$ per unit of
time. Let $\rho $ be the rate of interest. We will assume for simplicity
that the project output can be considered as a single unit, so that the output
price $p$ completely represents the revenue from the project. The firm's
trigger prices are $p_{H}$ and $p_{L}$, with $p_{H}>p_{L}$, such that if $%
p<p_{L}$ the project should be abandoned, if already undertaken, while if $%
p>p_{H}$ the project should be undertaken. If the firm has no investment
active in the project and assumes the price of its output will not fall back
down, the investment will be made if the price is greater than the full
cost: $p>w+\rho k$. If a firm has the investment already in place, and the
price suddenly drops to a lower level, the project will be dropped if $p<w$,
which is the variable cost. We are therefore in the presence of a
\textquotedblleft natural\textquotedblright\ area of inactivity between $w$
and $w+\rho k$, in which an active firm will not drop out and an idle firm
will not invest: it is obvious that the existence of this area depends
crucially on the existence of sunk costs $k$. What emerges from the
introduction of any kind of stochasticity in the price dynamics is
hysteresis, or the phenomenon by which the trigger prices $p_{L}$ and $p_{H}$
yield a larger area of inactivity, with the lower trigger price being below $%
w$ and the upper trigger price above the natural boundary $w+\rho k$. \\

In the standard \cite{Dixit1989} setup the prices are subject to White Gaussian Noise, and the model is studied by
appealing to option pricing arguments. We now examine a similar
framework using the same notation, with all the quantities defined so far
being scalar and non-stochastic, in which a firm's entry and exit decisions
are based on a market price that follows super-diffusive dynamics.\\

Consider the dynamics of the market price, defined by the following system:
$$\begin{cases}
dp_{t}=p_{t}[\mu dt+ \sigma dZ_{t}], \\ 
dZ_{t}=\gamma \tanh (\gamma Z_{t})dt+ dW _{t},\\
p_0 = p^0, \qquad Z_0 = 0,
\end{cases}
$$%
starting at $t=0$, with $\mu ,\sigma >0$, $W_{t}$ being White Gaussian noise and $\gamma = \sqrt{2\lambda}$ the DMPS risk parameter. The first equation has the same structure as the geometric Brownian motion, but the underlying dynamics are now driven by the DMPS noise source. The  will allow an explicit treatment of the firm's net value function $V_{i}(p)$, leading to
the conclusion that the system driven by a DMPS noise source has
entry and exit trigger prices that increase the effect of hysteresis, and
create an area of inaction that is larger than the area generated by a
system driven by WGN.
\\

 The decision problem of the firm consists of two state
variables, the price $p_{t}$ and a discrete variable that indicates whether
the firm is active (1) or not (0). In state $(p,0)$ the firm decides whether
to enter or to remain idle, and define $V_{0}(p)$ as the expected net
present value of starting with a price $p$ in the idle state and following
optimal policies. $V_{1}(p)$ is defined in an analogous manner. By Proposition 4, the asset equilibrium then follows the ODE:
%
%
%
%
%
$$
V_{0}(p_t)^{\prime \prime }\frac{\sigma ^{2}}{2}p_t^{2}+V_{0}(p_t)^{\prime }p_t(\mu
+\sigma\mathcal{B}\gamma )-\rho V_{0}(p_t)=0,
$$
and the exit condition for the active firm (including the flow of operating
profit) becomes:

$$
V_{1}(p_t)^{\prime \prime }\frac{\sigma ^{2}}{2}p_t^{2}+V_{1}(p_t)^{\prime }p_t(\mu
+\sigma \mathcal{B}\gamma)-\rho V_{1}(p)=w-p_t.
$$%
We search for a solution of the
form $V_{i}=A_ip^{
\beta ' }+B_ip^{\alpha '}$ with $i=0,1$ (both ODEs have the same homogeneous part), and the exponents  read:
%
%

\begin{eqnarray}
\beta ^{\prime } &=&\frac{1}{2}\left[ (1-m^{\prime })+\left( (1-m^{\prime})^2 +
 \frac{8\rho }{\sigma^2 }\right) ^{1/2}\right] >0, \\
\alpha ^{\prime } &=&\frac{1}{2}\left[ (1-m^{\prime })-\left( (1-m^{\prime})^
2+\frac{8\rho }{\sigma^2 }\right) ^{1/2}\right] <0,
\end{eqnarray}%
with $m^{\prime }:=2(\mu + \sigma\mathcal{B}\gamma )/\sigma ^{2}$. Note how the presence of DMPS in the $\sigma \gamma \mathcal{B}$ term creates an effect that is different from the effect of $\sigma$: the increase in the sensitivity of the price process to noise, which is also the variance of the Brownian motion, now has an effect which is the opposite of $\gamma$, the risk parameter for the state variable $p_t$. The two terms and their relative exponents are respectively the value of the options of entering and exiting: as before, setting $\gamma = 0$ one returns to the standard framework.     The two ODEs become:
\begin{eqnarray}
V_{0}(p_t) &=&A_{0}^{ \alpha ^{\prime }}+B_{0}^{ \beta ^{\prime
}},  \label{ODE1}  \\
V_{1}(p_t) &=&A_{1}^{ \alpha ^{\prime }}+B_{1}^{ \beta ^{\prime
}}+\frac{p_t}{\rho -\mu -\mathcal{B}\gamma }-\frac{w}{\rho }.
\label{ODE2}
\end{eqnarray}%
%
%
\begin{figure}
\centering
\begin{subfigure}[b]{0.45\textwidth}
\includegraphics[width=\textwidth]{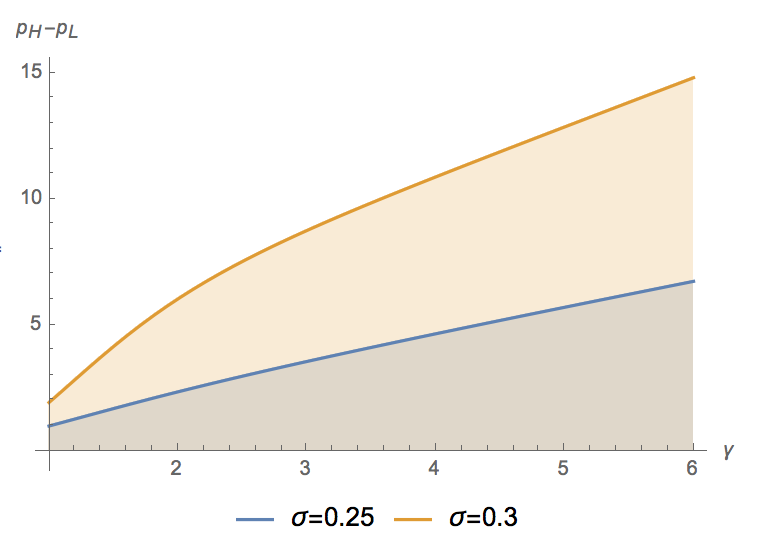}
\caption{An increase in risk}
\label{dixitfig}
\end{subfigure}
\begin{subfigure}[b]{0.45\textwidth}
\includegraphics[width=\textwidth]{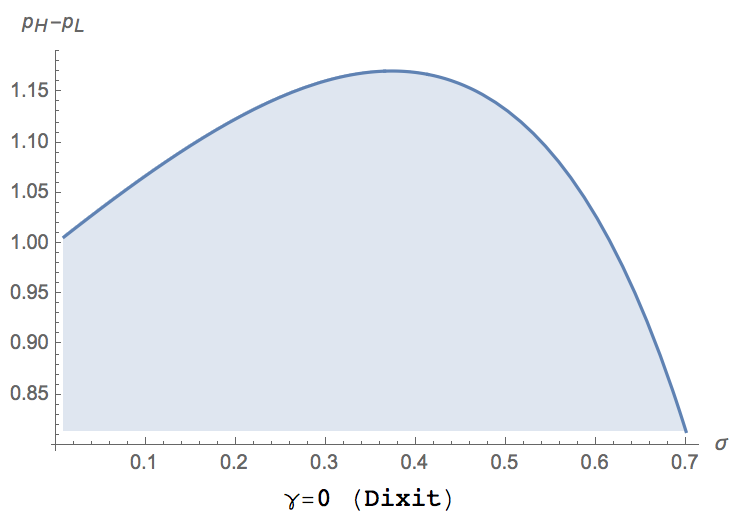}
\caption{An increase in variance}
\label{dixitorig}
\end{subfigure}
\caption{Impact of risk and variance on hysteresis}
\end{figure}
Note that the last two terms of \eqref{ODE2},$ \frac{p_t}{\rho -\mu - \sigma \mathcal{B}\gamma }-\frac{w}{\rho }$, if $\gamma=0$, are equal to: 
$$\mathbb{E}_t \int_0^\infty e^{-\rho s} (p_s - w) ds,$$ the discounted value at $t$ of keeping the project going until infinity. In our framework this value has a different interpretation: the firm cannot optimally decide until full information is achieved. It is  now the value of the project \emph{once the process $Z_t$ becomes Markovian}, which happens once the value of $\pm \gamma$ is realized. Intuitively, this depends on whether the drift randomly switches to positive or negative: a positive sign ($+\gamma$ realized) will increase the value for the firm to activate forever. Two endpoint conditions emerge naturally from the problem: $A_{0}=0$ and $%
B_{1}=0$ (if optimally idle, activating must be nearly useless, and vice
versa; we therefore write $A_{1}=A$ and $B_{0}=B$), allowing us to define
the two trigger prices $p_{H}$ and $p_{L}$ via two sets of conditions,
value-matching (equivalence in net present values) and smooth-pasting
(equivalence in the derivatives of the net present values). \ The first two
conditions are defined by:%
\begin{eqnarray}
A p_{L}^{\alpha ^{\prime }}+\frac{p_{L}}{\rho -\mu - 
\sigma\mathcal{B}\gamma }-\frac{w}{\rho } &=&B p_{L}^{\beta ^{\prime }}-l,
\label{valmatch1} \\
Ap_{H}^{\alpha ^{\prime }}+\frac{p_{H}}{\rho -\mu -
\sigma\mathcal{B}\gamma }-\frac{w}{\rho } &=&Bp_{H}^{\beta ^{\prime }}+k,
\label{valmatch2}
\end{eqnarray}%
and the smooth-pasting conditions are defined by:%
\begin{equation}
A\alpha ^{\prime }p_{i}^{\alpha ^{\prime }-1}+\frac{1}{\rho -\mu
-\sigma \mathcal{B}\gamma }-\frac{w}{\rho }=B\beta ^{\prime
}p_{i}^{\beta ^{\prime }-1},  \label{smoothpasting}
\end{equation}%
for $i=H,L$. With the four conditions (\ref{valmatch1})$-$(\ref%
{smoothpasting}) the pair of ODEs given by (\ref{ODE1}) and (\ref{ODE2}) is
completely determined. One analytical result is
worth noting: it is shown in \cite{Dixit1989} that $p_{H}>w+\rho
k\equiv W_{H}$ and $p_{L}<w-\rho l\equiv W_{L}$, which implies that
uncertainty increases the Marshallian area of inaction (full versus
relative costs), an interval of prices where an idle firm does not invest
and an active one does not exit. As $%
l $ grows, there is a finite price level that will result in the firm never exiting, which will be uniquely defined once the realization of $\mathcal{B}$ is observed, and $%
p_{H}$ will require $A=0$ in (\ref{valmatch2}) and (\ref{smoothpasting}). We
can then solve for $p_{H}$:

\begin{equation}
p_{H}^{\ast }=W_H\left( \frac{\beta ^{\prime }}{\beta ^{\prime }-1}%
\right) \left( \frac{\rho -\mu- \sigma\mathcal{B}\gamma }{\rho }\right) .\label{dixitph}
\end{equation}%
Similarly, if $k$ goes to infinity, the entry option becomes worthless and $%
B\rightarrow 0$. Solving for $p_{L}$, remembering that $\alpha'$ is negative:%

\begin{equation}
p_{L}^{\ast }=W_L \left( \frac{-\alpha ^{\prime }}{-\alpha ^{\prime }+1}%
\right) \left( \frac{\rho -\mu -\sigma \mathcal{B}\gamma }{\rho }\right).\label{pl}
\end{equation}%
Setting $\gamma=0$ and thereby reverting to the original Gaussian case considered by \cite{Dixit1989}, we recover his equations (23) and (24), p. 630. \\

Once we distinguish riskiness from variance by escaping the Gaussian setup, we can disentangle the two effects with a simple numerical simulation that shows that regardless of the sign of the realization of the Bernoulli $\pm \gamma$ the distance $p^*_H - p^*_L$ increases as $\gamma$ increases. The simulation shows that in the classical Dixit framework, obtained when $\gamma=0$, the relationship between the inaction area $p^*_H-p^*_L$ and $\sigma$ is concave and eventually slopes downwards, which is an unattractive (and surprisingly neglected) feature of the original model. In contrast, within a DMPS framework in which the risk of the price process is parameterized by $\gamma$, an increase in risk $d\gamma >0$  always causes an increase in inefficiency, as shown in Figure \ref{dixitfig}. For comparison purposes, we use the same parameter values used in the original Dixit formulation: $\mu = 0, \rho = 0.25, w=1, l=0,k=4$. 
Note that if $\sigma \to 0$ there is no uncertainty and therefore no hysteresis. It should also be noted that the \emph{joint} effect of increasing $\sigma$ and $\gamma$ causes a more than proportional expansion of the inaction area: in other words, even a little extra uncertainty has an even larger impact on hysteresis.


\section{Conclusions}

Rothschild and Stiglitz's concept of a Mean-Preserving Spread is the basic manner of rigorously characterizing an increase in risk in economics. In this paper, we extend this essential tool of economic theory to a dynamic setting. We define dynamic equivalents of the two original integral conditions in terms of transition cumulative densities and we prove a sufficiency condition for their satisfaction. We then focus on a specific class of non-Gaussian diffusion processes, the ballistic super-diffusive process. We prove the remarkable property, for a broad class of processes, that this process is unique in terms of satisfying the dynamic integral conditions. Moreover, this correlated noise source is shown to be the sum of two Gaussian processes with alternate drift: the result is a highly tractable analytical tool that allows one to escape the Gaussian straightjacket.  We characterize the probabilistic properties of three stochastic processes commonly used in economics (scalar, mean-reverting, geometric) driven by a DMPS process instead of the Gaussian, and provide a modified It\^o formula.  In our selection of economic illustrations which correspond to four canonical models widely used in the literature (portfolio choice, firm investment, asset dynamics and firm entry and exit), we show how moving beyond the Gaussian framework is both analytically tractable and intuitively important: this is because there is an unfortunate tendency in the profession to conflate risk with variance, which we disentangle thanks to our DMPS formulation.
%
%

\appendix

\section{Proof of Proposition 1}
\label{proofprop1}
$(i)$ By substituting (\ref{SCONDO}) into (\ref{MPSCONDO}), we can write:%
\begin{equation}
\qquad  {\frac{\partial }{\partial \lambda }}\left[ \int_{\mathbb{R}}%
\mathcal{P}^{(\lambda)}(x,t)dx\right] ={\frac{e^{-\lambda t}}{h^{(\lambda) }(0)%
}}\int_{\mathbb{R}}
\int_{-\infty }^{x}\left[ R^{(\lambda )}(y)q(0,0|y,t)\right] dy dx. 
\end{equation}
Similarly, considering \eqref{MPSCONDO2} and letting:

$$
\Psi^{(\lambda)} (x,t) = \int_{-\infty }^{x}\left[ R^{(\lambda )}(y)q(0,0|y,t)\right] dy,
$$
we can write: 
\begin{equation}
 \Phi ^{(\lambda)}(x,t):={\frac{\partial }{\partial \lambda }}\left[ \int_{-\infty }^{x}\mathcal{P}^{(\lambda)}(y,t)dy\right] ={\frac{e^{-\lambda t}}{h^{(\lambda)}(0)}}\int_{-\infty }^{x}\Psi ^{(\lambda )}(y,t)dy. \label{INTROD}
\end{equation}
From the condition given by equation (\ref{SCONDO}) and the fact that $%
q(0,0|x,t)=q(0,0|-x,t)$, we conclude that $\Psi ^{(\lambda )}(x,t)=-\Psi
^{(\lambda )}(-x,t)$ and its integral over $\mathbb{R}$ vanishes, leading to
the fulfillment of the first integral condition. \\ 

$(ii)$ Now consider the
curvature $\rho ^{(\lambda )}(x,t)$ of $\Phi ^{(\lambda )}(x,t)$, which reads:%
\begin{equation}
\rho ^{(\lambda )}(x,t)={\frac{\partial ^{2}}{\partial x^{2}}}\left[ \Phi
^{(\lambda )}(x,t)\right] =\left[ R^{(\lambda )}(x)q(0,0|x,t)\right] \geq 0.
\label{CURVOS}
\end{equation}%
From equation (\ref{INTROD}) we know that $\Phi ^{(\lambda
)}(\infty,t )=0$. This can only be achieved if we have $\Phi ^{(\lambda
)}(x,t)\geq 0$. The second integral condition is therefore verified.

\section{Proof of Lemma 1}
\label{prooflemma1}
Equation \eqref{DUALSDA} allows for the process $\hat{X}_t$ ---a solution of \eqref{SDA}, to admit a solution continuously in the entire time interval $[0,T]$; the problem is well-posed and therefore this condition is maintained with regularity for the subsequent $\lambda$-parameterization given by \eqref{DUALTPDA}. If a stochastic process $X_t^{(\lambda)}$ satisfies the sufficient condition \eqref{SCONDO}, then it also trivially satisfies the original Rothschild-Stiglitz integral conditions \eqref{MPSRS1} and \eqref{MPSRS2} if stopped at an arbitrary time $s \in [0,T]$. Therefore, for any $u(x,t)$ such that $u_{xx} \leq 0$ we have $u(X_t^{(\lambda_1)},s) > u(X_t^{(\lambda_2)},s)$ if $\lambda_1<\lambda_2$. Since equation \eqref{DUALTPDA} holds in the entire time domain $[0,T]$, an agent with time-consistent and time-invariant preferences has utility $u(x,t)$ such that for all $x_1, x_2, \alpha, \beta \in \mathbb{R}, \Delta_1, \Delta_2 \in \mathbb{R}^+, t, t' \in [0,T]$ the following holds, as in \cite{Halevy2015}:
\begin{eqnarray*}
u(x_1, t+\Delta_1) \sim_t u(x_2, t+\Delta_2) &\Leftrightarrow& u(x_1, t'+\Delta_1) \sim_{t'} u(x_2, t'+\Delta_2),\\
u(x_1, t+\Delta_1) \sim_t u(x_2, t+\Delta_2) &\Leftrightarrow& u(x_1, t+\Delta_1) \sim_{t'} u(x_2, t+\Delta_2),
\end{eqnarray*}
where $\sim_t$ means the symmetric ordering of the agent's temporal payments decided at time $t$. The Lemma then follows immediately.

\section{Proof of Proposition 2}
\label{proofprop2}

The first part of the proof is adapted from Theorem 1 in \cite{benjamini1997}. In the definition of \eqref{SDA}, we have assumed the drift $b(.)$ to have bounded and continuous first two derivatives: we can therefore condition the process $X_t$ such that $X_0 = X_T = a \in \mathbb{R}$. We can solve for $X_t$ in the following way: let $\hat{W}_t$ be a Brownian motion on the filtered probability space $(\Omega, \mathcal{F}, P)$ and $P_t$ the measure restricted to $\mathcal{F}_t$. Under $P_t$ we have that: 

$$
Z_t = \exp \left (\int_0^t b(\hat{W}_t) d\hat{W}_t - \frac{1}{2} \int_0^t b(\hat{W}_s)^2 ds \right ), 
$$
is a martingale by Girsanov's theorem. Define $B(.)$ as the antiderivative of $b(.)$. For the functional $B(t)$,  It\^o's formula implies $B(t) = B(0) + \int_0^t b(\hat{W}_t) d\hat{W}_t +  \frac{1}{2} \int_0^t b'(\hat{W}_s)ds$ and therefore:
$$
Z_t   = \exp \left ( B(t) - B(0) - \frac{1}{2} \int_0^t [ b'(\hat{W}_s) + b(\hat{W}_s)^2 ] ds \right ).
$$
We can see that $Z_t$ is the Radon-Nikodym derivative for the change of measure $dQ_t = Z_t dP_t$ such that under $Q$, the Brownian motion $\hat{W}_t$ is a solution of \eqref{SDA}. Now, if the process $X_t$ is a Brownian bridge, then $ b'(x) + b(x)^2$ must be constant, i.e. $ b'(W_s) + b(W_s)^2 = a \in \mathbb{R}$. A quick calculation shows that the only solution is $b(x) = a \tanh (ax + c)$ with $a, c$ constants. Setting $a = \sqrt{2\lambda}, c = 0$ immediately yields \eqref{prop2}. \\

It remains to prove that the stochastic process given by \eqref{prop2} satisfies the antisymmetry and positivity integral conditions by means of the sufficient condition given by \eqref{SCONDO}.  Let us now apply Proposition 1 to the pure Brownian motion obtained when $b(X_t)=0$ and $\sigma =1$ in
equation (\ref{SDA}). \ Accordingly, $\mathcal{L}_{x}(\cdot )={\frac{1}{2}}{%
\frac{\partial ^{2}}{\partial x^{2}}}(\cdot )$ and equation (\ref{STURM})
and its positive class of solutions read:%
\begin{equation}
{\frac{1}{2}}{\frac{\partial ^{2}}{\partial x^{2}}}h^{(\lambda )}(x)=\lambda
h^{(\lambda )}(x)\quad \Rightarrow \quad h^{(\lambda )}(x)=A\cosh \left( 
\sqrt{2\lambda }x\right) ,  \label{STURMA}
\end{equation}%
where $A$ is an arbitrary constant. Here we have $R^{(\lambda )}(x)=A\sinh
\left( \sqrt{2\lambda }x\right) $ and equation (\ref{SCONDO}) reads:%
\begin{equation}
R^{(\lambda )}(x)={\frac{d}{d\lambda }}A\cosh \left( \sqrt{2\lambda }%
x\right) ={\frac{Ax}{2\sqrt{2\lambda }}}\sinh \left[ \sqrt{2\lambda }x\right]
=R^{(\lambda )}(-x).  \label{SCONDOA}
\end{equation}%
In view of equation (\ref{DUALSDA}), we conclude that the resulting process
reads:%
\begin{equation}
dX_{t}=\left\{ \sqrt{2\lambda }\tanh \left[ \sqrt{2\lambda }X_{t}\right]
\right\} dt+dW_{t},  \label{TANHL}
\end{equation}%
which proves Proposition 2.

\section{Proof of Lemma 2}

Writing the transformation:

$$
P(x,t | x_0, 0) = e^{-\lambda t} \cosh (\sqrt{2\lambda} x) Q(x,t | x_0, 0),
$$
one can see that equation \eqref{kolm} reduces to the Kolmogorov Forward equation for a standard Brownian motion:

$$
\frac{\partial}{\partial t} P(x,t|x_0,0) = \frac{1}{2} \frac{\partial^2}{\partial x^2} P(x,t| x_0,0).
$$ 
It is well known that this linear PDE (the heat equation) has the Gaussian solution:
$$
P(x,t|x_0,0) = \frac{1}{\sqrt{2\pi t}} e^{-\frac{(x-x_0)^2}{2t}}.
$$
Reverting to the measure $Q(.)$, and noticing that $\cosh (x) = \frac{e^x + e^{-x}}{2}$, after rearranging one immediately obtains: 

$$
Q(x,t|x_0,0) = \frac{1}{2\sqrt{2\pi t}} \left (e^{-\frac{(x+\sqrt{2\lambda} t)^2}{2t}} + e^{-\frac{(x-\sqrt{2\lambda} t)^2}{2t}} \right ), 
$$
which proves \eqref{TPDTANHL}. This equation shows how the density of the process \eqref{prop2} is the average of two $\pm \sqrt{2\lambda}$-drifted Gaussian densities with unit variance. But this can be read as the expected value of a Bernoulli variable: the process has 0.5 probability of having a density with positive drift and 0.5 probability of having the negative one, as first noted by \cite{Rogers1981}.\footnote{%
See \textit{their} example 2.} The process \eqref{prop2} can therefore be rewritten as: 

$$
dX_t = \mathcal{B} dt + dW_t,
$$
with $\mathcal{B}$ a Bernoulli variable taking values $\pm \sqrt{2\lambda}$ with probability 0.5, which proves Lemma 2.

\section{Proof of Proposition 3}

\label{pro4}
\underline{Proof of the drifted scalar case}.  Writing: 
$$
P(z,x,t\mid z_{0},0,0)=e^{-\lambda t}\cosh \left[ \sqrt{2\lambda }x%
\right] Q(z,x,t\mid y_{0},0,0),
$$ 
one immediately verifies that with the rescaling $\sigma \hat{z} = z - \mu t$, the Kolmogorov equation \eqref{FP2} reduces to a pure diffusion equation on $\mathbb{R}\times \mathbb{R}$:
\begin{equation}
\begin{array}{l}
\partial _{t}Q( \hat{z},x,t\mid y_{0},0,0)=\left\{ {\frac{1}{2}}\partial _{\hat{z}%
 \hat{z}}+\partial _{ \hat{z}x}+{\frac{1}{2}}\partial _{xx}\right\} Q( \hat{z}%
,x,t\mid z_{0},0,0).%
\end{array}
\label{FPBSDYN1}
\end{equation}%
This equation identifies a Gaussian bivariate TPD with canonical structure:%

\begin{equation}
\left\{ 
\begin{array}{l}
Q( \hat{z},x,t\mid z_{0},0,0)={\frac{1}{2\pi \sqrt{\Delta (t)}}}e^{-{\frac{1}{%
2\Delta (t)}}\left[ a(t) \hat{z}^{2}-2h(t) \hat{z}x+b(t)x^{2}\right] }, \\ 
\Delta (t)=a(t)b(t)-h(t)^{2}.%
\end{array}%
\right.  \label{Q}
\end{equation}%
The marginal TPD $P_{M}(z,t\mid z_{0},0)$ is obtained from the following quadrature:

\begin{eqnarray}
P_{M}( \hat{z},t\mid z_{0},0)&=&\int_{\mathbb{R}}e^{-\lambda t}\cosh \left[ \sqrt{%
2\lambda }\ x\right] Q( \hat{z},x,t\mid \hat{z}_{0},0,0)dx, \nonumber \\ 
 &=&{\frac{e^{-{\frac{a(t)}{2\Delta (t)}%
} \hat{z}^{2}}e^{-\lambda t}}{2\pi \sqrt{\Delta (t)}}}\left\{ \sqrt{{\frac{%
2\pi \Delta (t)}{b(t)}}}e^{\frac{h(t)^{2} \hat{z}^{2}}{2\Delta (t)b(t)}}\cosh %
\left[ \sqrt{2\lambda }{\frac{h(t)}{b(t)}} \hat{z}\right] e^{{\lambda {\frac{%
\Delta (t)}{b(t)}}}}\right\}, \nonumber \\ 
& =&{\frac{1}{\sqrt{2\pi b(t)}}}e^{-{%
\frac{ \hat{z}^{2}}{2b(t)}}}\cosh \left[ \sqrt{2\lambda }{\frac{h(t)}{b(t)}}%
 \hat{z}\right] e^{\lambda \left[ {\frac{\Delta (t)}{b(t)}}-t\right] }.\label{PMTPD}
\end{eqnarray}
Now, in our case we have:%
\begin{equation}
a(t)=b(t)=h(t)=t\quad \Rightarrow \quad \Delta (t)=a(t)b(t)-h^{2}(t)=0,
\label{CHECK}
\end{equation}%
and therefore equation (\ref{PMTPD}) can be rewritten as:%
\begin{equation}
P_{M}( \hat{z},t\mid y_{0},0)={\frac{1}{\sqrt{2\pi t}}}\left\{ e^{-{\frac{( \hat{z}-
\sqrt{2\lambda }t)^{2}}{2t}}}+e^{-{\frac{( \hat{z}+\sqrt{2\lambda }t)^{2}}{2t}
}}\right\} .  \label{COSSUM}
\end{equation}%
The coefficients 
$a(t)$, $b(t)$ and $h(t)$ given in equation (\ref{CHECK}) follow by using
the \cite{Chandrasekhar1943} general procedure.\footnote{%
Here one directly uses Lemma II and equations (260)-(263) from \cite{Chandrasekhar1943}.} Reverting from the rescaling to the original variable $Z_t$ one immediately obtains \eqref{prop3}, and Proposition 3.1 is proven. \\  

Let us remark
from equation (\ref{Q}) that $\Delta (t)=0$ indicates that the $(\hat{Z}%
_{t},X_{t})$ is actually a degenerate diffusion process in $\mathbb{R}^{2}$.
This can be understood by observing that the stochastic differential
equation corresponding to the diffusion operator (\ref{FPBSDYN1}) reads:%
\begin{equation}
\left\{ 
\begin{array}{l}
d\hat{Z}_{t}=dW_{t}, \\ 
dX_{t}=dW_{t}.%
\end{array}%
\right.  \label{DEGENERATE}
\end{equation}%
This shows that $\hat{Z}_{t}=X_{t}+\mathrm{const}=W_{t}$, and the underlying
dynamics degenerate to a scalar Wiener process on the $\mathbb{R}^{2}$
plane. \\

\underline{Proof of the stationary Ornstein-Uhlenbeck case}. Using the Bernoulli Representation Lemma, one can solve the KFE separately for each of the $\pm \sqrt{2\lambda}$ realizations. Writing \eqref{KOLMOOP} for each $\mu(z) = \alpha(\mu \pm \sqrt{2\lambda} - z)$, one can solve the equation by taking the Fourier transform in $z$ and then take the inverse transform; one immediately obtains each of the $\mu \pm \sqrt{2\lambda}$-reverting O-U densities. This procedure (one for each Bernoulli realization) is lengthy but standard, and is therefore omitted. \\

\underline{Proof of the geometric case}.  Defining $Y_{t}=\log (Z_t)$ the process $\{ Z_{t},X_{t} \}$ is a
diffusion process on $\mathbb{R}^{+}\times \mathbb{R}$ and the associated
Kolmogorov equation (the underlying stochastic integrals are interpreted
in the It\^{o} sense) reads: 
\begin{equation}
\partial _{t}P(z,x,t\mid z_{0},0,0) = \mathcal{F} (P(z,x,t\mid y_{0},0,0)),
\label{FPBSDYN}
\end{equation}
where the operator $\mathcal{F}(.)$ is given by \eqref{KOLMOOP}. One proceeds as in the previous proof, then reverts to the original scaling and one obtains immediately \eqref{prop3geom} and Proposition 3.3 is proved.\\

\section{Proof of Proposition 4}
\label{itolemma}
Since $X_t$ is an It\^o process,  by means of the Bernoulli Representation Lemma, equation \eqref{itodmps} is a straightforward application of It\^o's Lemma. The more general formula for a function $g(Z_t, X_t, t) \in \mathcal{C}^{2,2,1}$ that is explicitly a function of the noise is given by:
\begin{eqnarray}
d g(Z_t, X_t,t) &=  & \left [ \frac{\partial g}{\partial t} + \bigg ( \mu(Z_t, t) +  \sigma(X_t, t) \sqrt{2\lambda } \tanh(\sqrt{2\lambda }X_t ) \bigg) \frac{\partial g}{\partial z} + \right.  \nonumber \\
 & +&\left .  \sqrt{2\lambda }\tanh(\sqrt{2\lambda }X_t)\frac{\partial g}{\partial x}+ \frac{\sigma(Z_t,t)^2}{2}\frac{\partial^2 g}{\partial z^2}  + \sigma(Z_t,t) \frac{\partial^2 g}{\partial z \partial x} + \frac{1}{2} \frac{\partial^2 g}{\partial x^2}  \right ]dt  \nonumber \\
&+&  \left (\sigma(Z_t,t) \frac{\partial g}{\partial z} + \frac{\partial g}{\partial x} \right ) dW_t,\end{eqnarray}%
by application of the multidimensional It\^o formula.

\section{Proof of Proposition 5}
\label{proofprop5}
When $f(z)=-f(-z)$, leading to $%
P_{sM}(z)=P_{sM}(-z)$, the curvature at the origin $\rho
_{(\sigma ,\lambda )}=\left[ {\frac{\partial ^{2}}{\partial x^{2}}}P_{sM}(x)%
\right] \mid _{x=0}$ is given by:
\begin{eqnarray}
\rho _{(\sigma ,\lambda )}(z) &=&  \mathcal{N} \frac{\partial^2}{\partial z^2} \exp \left [
\frac{2F(z) + \log \cosh (\sqrt{2\lambda}z)}{\sigma ^{2}} \right ],  \nonumber  \\
&=& \left \{ \left ( \frac{2 f(z)}{\sigma^2} + \sqrt{2\lambda z}  \tanh (\sqrt{2\lambda}z )  \right )^2 + \left ( \frac{2 f'(z)}{\sigma^2}+ \frac{2 \lambda}{\cosh(\sqrt{2 \lambda}z)} \right )   \right \}  P_s (z),\nonumber \\
\rho _{(\sigma ,\lambda )}(0) &=& \underset{:=\rho _{(\sigma ,0)}}{\underbrace{{%
\frac{2}{\sigma ^{2}}}\left[ {\frac{\partial }{\partial z}}f(z)\right] \mid
_{z=0}}}+2\lambda =\rho _{(\sigma ,0)}(0) +2\lambda ,  \label{CURVA}
\end{eqnarray}
where $\rho _{(\sigma ,0)}(0)$ is the curvature at the origin of the Gaussian measure. This proves Proposition 5.

\bibliographystyle{chicago}

\bibliography{my_refs.bib}

\end{document}